\documentclass[reqno]{amsart}

\usepackage{amsmath,amsfonts,txfonts,amscd,amsthm,graphicx,amssymb,epic,eepic}

 \oddsidemargin 25pt \evensidemargin 25pt

\theoremstyle{plain}
\newtheorem{theorem}{Theorem}
\newtheorem{proposition}{Proposition}
\newtheorem{corollary}{Corollary}
\newtheorem{lemma}{Lemma}
\numberwithin{equation}{section}

\newcommand{\eps}{\varepsilon}

\newcommand{\TT}{\mathcal T}

\newcommand{\Z}{\mathbb Z}
\newcommand{\FF}{\mathcal F}
\newcommand{\N}{\mathbb N}
\newcommand{\R}{\mathbb R}

\renewcommand{\P}{\mathbb P}

\newcommand{\EE}{\mathcal E}

\newcommand{\II}{\mathcal I}
\newcommand{\JJ}{\mathcal J}

\newcommand{\Fl}{\FF^{(\ell)}}

\title[The linear flow in a honeycomb]
{On the distribution of the free path length of \\ the linear flow
in a honeycomb}

\author[F.\,P. Boca and R.\,N. Gologan]{Florin P. Boca and Radu N. Gologan}

\date{June 7, 2008}

\address{Department of Mathematics, University of Illinois at
Urbana-Champaign, 1409 W. Green St., Urbana, IL 61801, USA}

\address{Institute of Mathematics of the Romanian Academy, P.O.Box 1-764,
Bucharest RO-014700, Romania}

\email{fboca@math.uiuc.edu}

\address{Institute of Mathematics of the Romanian Academy, P.O. Box 1-764,
Bucharest RO-014700, Romania}

\email{Radu.Gologan@imar.ro}

\begin{document}

\begin{abstract}
Let $\ell \geqslant 2$ be an integer. For each $\eps \in (0,\frac{1}{2})$
remove from $\R^2$ the union of discs of radius $\eps $ centered
at the integer lattice points $(m,n)$, with $m\nequiv n \pmod{\ell}$.
Consider a point-like particle moving
linearly at unit speed, with velocity $\omega$, along a trajectory
starting at the origin, and its free path length $\tau_{\ell,\eps}
(\omega)\in [0,\infty]$. We prove the weak convergence of the
probability measures associated with the random variables $\eps
\tau_{\ell,\eps}$ as $\eps \rightarrow 0^+$ and explicitly compute
the limiting distribution. For $\ell=3$ this leads to an
asymptotic formula for the length of the trajectory of a billiard
in a regular hexagon, starting at the center, with circular
pockets of radius $\eps\rightarrow 0^+$ removed from the corners.
For $\ell=2$ this corresponds to the trajectory of a billiard in a
unit square with circular pockets removed from the corners and
trajectory starting at the center of the square. The limiting
probability measures on $[0,\infty)$ have a tail at infinity,
which contrasts with the case of a square with pockets and
trajectory starting from one of the corners, where the limiting
probability measure has compact support.
\end{abstract}

\maketitle

\section{Introduction}
Recent progress led to a better understanding of the statistics of
the free path length of the periodic Lorentz gas in the small
scatterer limit \cite{BGZ1,BGZ2,BZ2,CG1,CG2,Dah,G,MS}. The aim of this
paper is to study situations where periodicity conditions are
altered by imposing certain congruence conditions on the integer
lattice points where scatterers are placed. The case of the
honeycomb lattice arises as a particular example in this context.
Here we only consider the situation where motion originates at the
origin, extending the results from \cite{BGZ1} and \cite{BGZ2}.
The case where the initial position is randomly chosen is more
intricate and will not be treated here.

Let $\ell \geqslant 2$ be an integer. Consider the set
$\Z^2_{(\ell)}$ of pairs of integers $(m,n)$ with $m\nequiv
n \pmod{\ell}$. For every $\eps >0$ consider the ``fat
lattice points" (\emph{scatterers}) given by small discs of
radius $\eps$ centered at all points of $\Z^2_{(\ell)}$ and the
region
\begin{equation*}
Z_{\ell,\eps}=\{ x\in\R^2 :\operatorname{dist}
(x,\Z^2_{(\ell)})\geqslant \eps \}
\end{equation*}
obtained by removing all scatterers. In $\R^2$ consider a
point-like particle moving at constant unit speed along a linear
trajectory originating at $(0,0)$. The \emph{free path length}
(\emph{first exit time}) is defined as
\begin{equation*}
\tau_{\ell,\eps} (\omega)=\inf \{ \tau >0: \tau \vec{\omega} \in
\partial Z_{\ell,\eps} \},
\end{equation*}
the distance traveled to reach the first scatterer along the
direction $\vec{\omega}=e^{i\omega}\in {\mathbb T}$, $\omega\in [0,2\pi)$, and as
$+\infty$ when the particle escapes to infinity without reaching
any scatterer. The Lebesgue measure of a measurable set $A\subseteq \R$ is denoted by
$\vert A\vert$. This paper is concerned with the study, in the
small scatterer limit ($\eps\rightarrow 0^+$), of the asymptotic
behavior of the repartition function of $\eps \tau_{\ell,\eps}$
defined by
\begin{equation*}
{\mathbb P}_{\ell,\eps} (\lambda)=\frac{1}{2\pi} \bigg|\bigg\{
\omega \in [0,2\pi):
\tau_{\ell,\eps}(\omega)>\frac{\lambda}{\eps}\bigg\} \bigg|.
\end{equation*}

To accomplish this we first consider the situation where
scatterers are obtained by translating the vertical segment
$V_\eps =\{0\}\times [-\eps,\eps]$ by $(m,n)\in\Z^2_{(\ell)}$ and
estimate, for any interval $I\subseteq [0,1]$ of length $\vert
I\vert \asymp \eps^{c}$ with fixed $c\in (0,1)$, the repartition
\begin{equation*}
{\mathbb G}_{\ell,I,\eps}(\lambda)=\left| \left\{ \omega\in
\arctan I: q_{\ell,\eps} (\omega) >\frac{\lambda}{\eps} \right\}
\right|,\qquad \eps\rightarrow 0^+,
\end{equation*}
of the \emph{horizontal free path length}
\begin{equation*}
q_{\ell,\eps} (\omega)=\inf\left\{ q: (q,q\tan\omega)
\in \Z^2_{(\ell)}+V_\eps \right\},\qquad \omega \in \bigg[
0,\frac{\pi}{4} \bigg].
\end{equation*}

In this paper $\varphi$ will denote Euler's totient function.
The dilogarithm is defined by
\begin{equation*}
\operatorname{Li}_2 (x)=\sum\limits_{n=1}^\infty \frac{x^n}{n^2} =
-\int_0^x \frac{\ln (1-t)}{t}\, dt,\qquad x\in [0,1].
\end{equation*}
Clearly $\operatorname{Li}_2 (1)=\zeta (2)=\frac{\pi^2}{6}$.

The main result of this paper shows that the limit of $\P_{\ell,\eps}$
exists as $\eps\rightarrow 0^+$ and this limit is explicitly computed.

\begin{theorem}\label{T1}
{\em (i)} For every $0<c_1<1$ and $\delta >0$, as $\eps
\rightarrow 0^+$,
\begin{equation}\label{1.1}
{\mathbb G}_{\ell,I,\eps} (\lambda)=c_I G_\ell (\lambda) +O_{\delta,\lambda,\ell}
\left( \eps^{-\delta+\theta (c,c_1)}\right),\qquad \lambda >0,
\end{equation}
where
\begin{equation*}
c_I=\int_I \frac{du}{1+u^2},\qquad \theta (c,c_1)=\min \bigg\{
c+c_1,\frac{1}{2}-2c_1\bigg\},
\end{equation*}
and the limiting repartition function $G_\ell$ is given by
\begin{equation*}
G_\ell (\lambda) =\begin{cases} \vspace{.1cm} 1-\left(
\frac{1}{\zeta (2)}+A(\ell)\right) \lambda & \mbox{\rm
if $\ \lambda \in \big( 0, \frac{1}{2}\big],$} \\
\vspace{.1cm} 1-\frac{\lambda}{\zeta (2)}+A(\ell) H_2 (\lambda) &
\mbox{\rm if
$\ \lambda\in\big[ \frac{1}{2}, 1\big],$} \\
\frac{2C(\ell)}{\ell} H_3 (\lambda) & \mbox{\rm if $\
\lambda\in [1,\infty),$}
\end{cases}
\end{equation*}
with
\begin{equation}\label{1.2}
\begin{split}
C(\ell) & =\frac{\varphi (\ell)}{\zeta (2)\ell}
\prod\limits_{\substack{p\mid \ell \\ p \operatorname{prime}}}
\left( 1-\frac{1}{p^2}\right)^{-1}=\frac{\varphi (\ell)}{\ell}
\prod\limits_{\substack{p\nmid \ell \\ p \operatorname{prime}}}
\left( 1-\frac{1}{p^2}\right) ,
\quad A(\ell)=\frac{1}{\zeta (2)}-\frac{2C(\ell)}{\ell}, \\
H_2 (\lambda) & =3\lambda-2+\zeta (2)-(\ln \lambda)^2
+2(1-\lambda)\ln \left(\frac{1}{\lambda}-1\right)-2
\operatorname{Li}_2 (\lambda) ,\\ H_3 (\lambda) &
=\operatorname{Li}_2 \left( \frac{1}{\lambda}\right)
-(\lambda-1)\ln \left( 1-\frac{1}{\lambda}\right) -1.
\end{split}
\end{equation}

{\em (ii)} For every $\delta >0$, as $\eps\rightarrow 0^+$,
\begin{equation*}
{\mathbb P}_{\ell,\eps} (\lambda)=G_\ell (\lambda)+O_{\delta,\lambda,\ell}
\left(\eps^{\frac{1}{8}+\delta}\right),\qquad \lambda >0.
\end{equation*}
\end{theorem}

\begin{figure}[ht]
\includegraphics*[scale=0.65, bb=0 0 280 165]{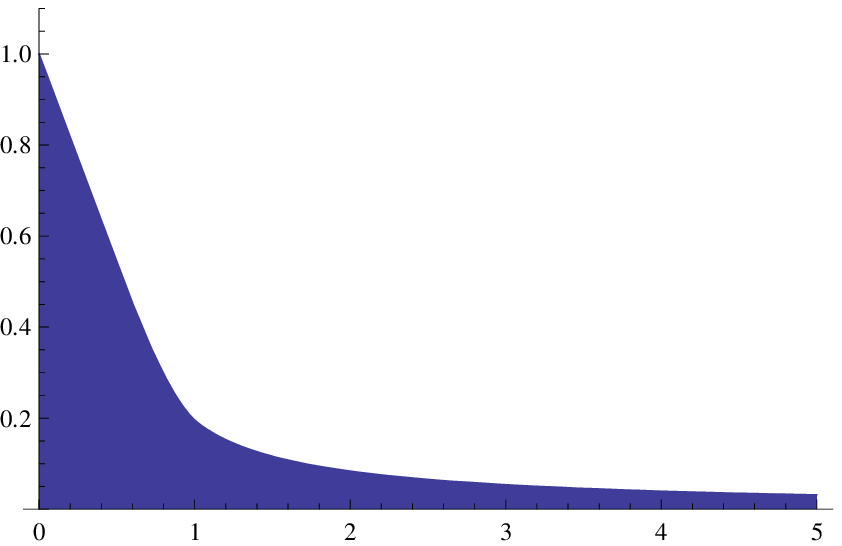}
\includegraphics*[scale=0.65, bb=-20 0 250 165]{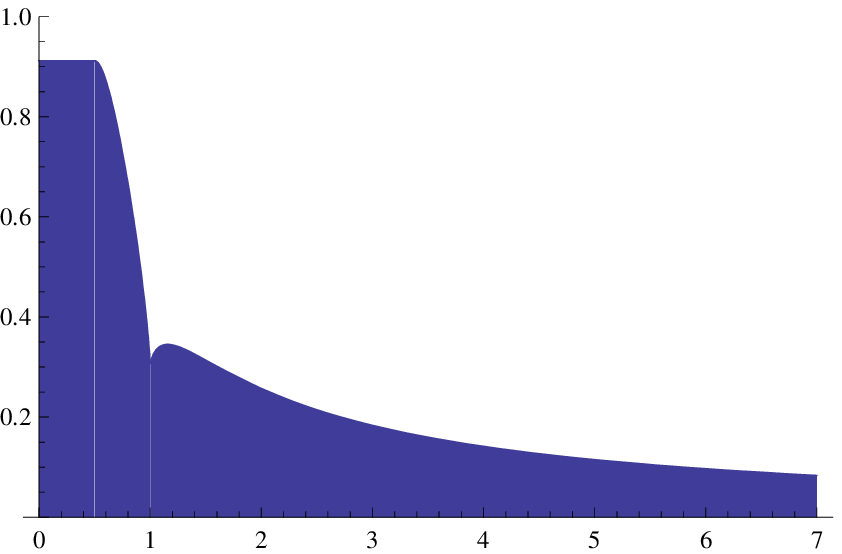}
\caption{\it \small The repartition function $G_3$ and the density
function $g_3$ } \label{Figure1}
\end{figure}

The continuity of $G_\ell$ at $\lambda=\frac{1}{2}$ is equivalent
with the well known dilogarithm identity
\begin{equation*}
\sum_{n=1}^\infty \frac{1}{2^n n^2}=\operatorname{Li}_2
\left(\frac{1}{2}\right)=\int_0^{1/2} \frac{1}{u}\
\ln\frac{1}{1-u}\ du =\frac{\zeta (2)-(\ln 2)^2}{2} .
\end{equation*}

Since $\frac{C(\ell)}{\ell}\rightarrow 0$ and $A(\ell)\rightarrow
\frac{1}{\zeta (2)}$ as $\ell\rightarrow \infty$, the compactly
supported limiting repartition $H(\lambda)$ from \cite[Theorem
1.1]{BGZ1} is being recovered as $\lim_{\ell\rightarrow\infty}
G_\ell (\lambda)$.

\begin{figure}[ht]
\begin{center}
\unitlength 0.25mm
\begin{picture}(550,270)(0,-43)

\dottedline{2}(90.43,16.576)(550,205)

\path(0,0)(25,-43.3)(75,-43.3)(100,0)(150,0)(175,-43.3)(225,-43.3)(250,0)(300,0)(325,-43.3)(375,-43.3)(400,0)(450,0)(475,-43.3)(525,-43.3)(550,0)
\path(0,0)(25,43.3)(75,43.3)(100,0) \path(150,0)(175,43.3)(225,43.3)(250,0)
\path(300,0)(325,43.3)(375,43.3)(400,0) \path(450,0)(475,43.3)(525,43.3)(550,0)

\path(0,86.6)(25,43.3) \path(75,43.3)(100,86.6)(150,86.6)(175,43.3)
\path(225,43.3)(250,86.6)(300,86.6)(325,43.3) \path(375,43.3)(400,86.6)(450,86.6)(475,43.3) \path(525,43.3)(550,86.6)

\path(0,86.6)(25,129.9)(75,129.9)(100,86.6) \path(150,86.6)(175,129.9)(225,129.9)(250,86.6)
\path(300,86.6)(325,129.9)(375,129.9)(400,86.6) \path(450,86.6)(475,129.9)(525,129.9)(550,86.6)

\path(0,173.2)(25,129.9) \path(75,129.9)(100,173.2)(150,173.2)(175,129.9)
\path(225,129.9)(250,173.2)(300,173.2)(325,129.9) \path(375,129.9)(400,173.2)(450,173.2)(475,129.9) \path(525,129.9)(550,173.2)

\path(0,173.2)(25,216.5)(75,216.5)(100,173.2) \path(150,173.2)(175,216.5)(225,216.5)(250,173.2)
\path(300,173.2)(325,216.5)(375,216.5)(400,173.2) \path(450,173.2)(475,216.5)(525,216.5)(550,173.2)

\put(90.43,16.576){\makebox(0,0){{\tiny $\circ$}}}
\put(95,25){\makebox(0,0){{\tiny $C_1$}}}
\put(171.29,49.73){\makebox(0,0){{\tiny $\circ$}}}
\put(175,58){\makebox(0,0){{\tiny $C_2$}}}
\put(21.29,-36.87){\makebox(0,0){{\tiny $\circ$}}}
\put(15,-39){\makebox(0,0){{\tiny $C_2$}}}

\put(246.52,80.573){\makebox(0,0){{\tiny $\circ$}}}
\put(252,73){\makebox(0,0){{\tiny $C_3$}}}
\put(96.52,-5.1){\makebox(0,0){{\tiny $\circ$}}}
\put(102,-12){\makebox(0,0){{\tiny $C_3$}}}
\put(261.22,86.6){\makebox(0,0){{\tiny $\circ$}}}
\put(262,94){\makebox(0,0){{\tiny $C_4$}}}
\put(94.39,9.717){\makebox(0,0){{\tiny $\circ$}}}
\put(103,10){\makebox(0,0){{\tiny $C_4$}}}

\put(312.03,107.43){\makebox(0,0){{\tiny $\circ$}}}
\put(318,101){\makebox(0,0){{\tiny $C_5$}}}
\put(50.94,43.3){\makebox(0,0){{\tiny $\circ$}}}
\put(51,52){\makebox(0,0){{\tiny $C_5$}}}

\put(366.83,129.9){\makebox(0,0){{\tiny $\circ$}}}
\put(366,122){\makebox(0,0){{\tiny $C_6$}}}
\put(4.1,7.074){\makebox(0,0){{\tiny $\circ$}}}
\put(-4,9){\makebox(0,0){{\tiny $C_6$}}}
\put(377.53,134.29){\makebox(0,0){{\tiny $\circ$}}}
\put(374.5,142){\makebox(0,0){{\tiny $C_7$}}}
\put(2.53,-4.39){\makebox(0,0){{\tiny $\circ$}}}
\put(-0.5,-9.5){\makebox(0,0){{\tiny $C_7$}}}

\put(454.3,165.76){\makebox(0,0){{\tiny $\circ$}}}
\put(453,158){\makebox(0,0){{\tiny $C_8$}}}
\put(535.15,198.91){\makebox(0,0){{\tiny $\circ$}}}
\put(534,191){\makebox(0,0){{\tiny $C_9$}}}
\put(79.3,-35.86){\makebox(0,0){{\tiny $\circ$}}}
\put(87,-39){\makebox(0,0){{\tiny $C_8$}}}
\put(10.15,17.59){\makebox(0,0){{\tiny $\circ$}}}
\put(2,21){\makebox(0,0){{\tiny $C_9$}}}
\put(50,0){\makebox(0,0){{\tiny $\star$}}}

\put(25,-43.3){\makebox(0,0){{\footnotesize $\bullet$}}}
\put(75,-43.3){\makebox(0,0){{\footnotesize $\bullet$}}}
\put(175,-43.3){\makebox(0,0){{\footnotesize $\bullet$}}}
\put(225,-43.3){\makebox(0,0){{\footnotesize $\bullet$}}}
\put(325,-43.3){\makebox(0,0){{\footnotesize $\bullet$}}}
\put(375,-43.3){\makebox(0,0){{\footnotesize $\bullet$}}}
\put(475,-43.3){\makebox(0,0){{\footnotesize $\bullet$}}}
\put(525,-43.3){\makebox(0,0){{\footnotesize $\bullet$}}}

\put(0,0){\makebox(0,0){{\footnotesize $\bullet$}}}
\put(100,0){\makebox(0,0){{\footnotesize $\bullet$}}}
\put(150,0){\makebox(0,0){{\footnotesize $\bullet$}}}
\put(250,0){\makebox(0,0){{\footnotesize $\bullet$}}}
\put(300,0){\makebox(0,0){{\footnotesize $\bullet$}}}
\put(400,0){\makebox(0,0){{\footnotesize $\bullet$}}}
\put(450,0){\makebox(0,0){{\footnotesize $\bullet$}}}
\put(550,0){\makebox(0,0){{\footnotesize $\bullet$}}}

\put(25,43.3){\makebox(0,0){{\footnotesize $\bullet$}}}
\put(75,43.3){\makebox(0,0){{\footnotesize $\bullet$}}}
\put(175,43.3){\makebox(0,0){{\footnotesize $\bullet$}}}
\put(225,43.3){\makebox(0,0){{\footnotesize $\bullet$}}}
\put(325,43.3){\makebox(0,0){{\footnotesize $\bullet$}}}
\put(375,43.3){\makebox(0,0){{\footnotesize $\bullet$}}}
\put(475,43.3){\makebox(0,0){{\footnotesize $\bullet$}}}
\put(525,43.3){\makebox(0,0){{\footnotesize $\bullet$}}}

\put(0,86.6){\makebox(0,0){{\footnotesize $\bullet$}}}
\put(100,86.6){\makebox(0,0){{\footnotesize $\bullet$}}}
\put(150,86.6){\makebox(0,0){{\footnotesize $\bullet$}}}
\put(250,86.6){\makebox(0,0){{\footnotesize $\bullet$}}}
\put(300,86.6){\makebox(0,0){{\footnotesize $\bullet$}}}
\put(400,86.6){\makebox(0,0){{\footnotesize $\bullet$}}}
\put(450,86.6){\makebox(0,0){{\footnotesize $\bullet$}}}
\put(550,86.6){\makebox(0,0){{\footnotesize $\bullet$}}}

\put(25,129.9){\makebox(0,0){{\footnotesize $\bullet$}}}
\put(75,129.9){\makebox(0,0){{\footnotesize $\bullet$}}}
\put(175,129.9){\makebox(0,0){{\footnotesize $\bullet$}}}
\put(225,129.9){\makebox(0,0){{\footnotesize $\bullet$}}}
\put(325,129.9){\makebox(0,0){{\footnotesize $\bullet$}}}
\put(375,129.9){\makebox(0,0){{\footnotesize $\bullet$}}}
\put(475,129.9){\makebox(0,0){{\footnotesize $\bullet$}}}
\put(525,129.9){\makebox(0,0){{\footnotesize $\bullet$}}}

\put(0,173.2){\makebox(0,0){{\footnotesize $\bullet$}}}
\put(100,173.2){\makebox(0,0){{\footnotesize $\bullet$}}}
\put(150,173.2){\makebox(0,0){{\footnotesize $\bullet$}}}
\put(250,173.2){\makebox(0,0){{\footnotesize $\bullet$}}}
\put(300,173.2){\makebox(0,0){{\footnotesize $\bullet$}}}
\put(400,173.2){\makebox(0,0){{\footnotesize $\bullet$}}}
\put(450,173.2){\makebox(0,0){{\footnotesize $\bullet$}}}
\put(550,173.2){\makebox(0,0){{\footnotesize $\bullet$}}}

\put(25,216.5){\makebox(0,0){{\footnotesize $\bullet$}}}
\put(75,216.5){\makebox(0,0){{\footnotesize $\bullet$}}}
\put(175,216.5){\makebox(0,0){{\footnotesize $\bullet$}}}
\put(225,216.5){\makebox(0,0){{\footnotesize $\bullet$}}}
\put(325,216.5){\makebox(0,0){{\footnotesize $\bullet$}}}
\put(375,216.5){\makebox(0,0){{\footnotesize $\bullet$}}}
\put(475,216.5){\makebox(0,0){{\footnotesize $\bullet$}}}
\put(525,216.5){\makebox(0,0){{\footnotesize $\bullet$}}}

\put(67,3.5){\makebox(0,0){\tiny $\omega$}}
\put(50,0){\arc{20}{-0.34}{0}}
\put(50,0){\arc{18}{-0.34}{0}}

\put(43,0){\makebox(0,0){\scriptsize $O$}}

\thinlines
\path(50,0)(90.43,16.576)(21.29,-36.87)(96.52,-5.1)(94.39,9.717)(50.94,43.3)(4.1,7.074)(2.53,-4.39)(79.3,-35.86)(10.15,17.59)(50,5.5)
\dottedline{2}(100,0)(50,0)



\thicklines \path(0,0)(25,43.3)(75,43.3)(100,0)(75,-43.3)(25,-43,3)(0,0)

\end{picture}
\end{center}
\caption{\it \small The free path in a hexagonal billiard and in a honeycomb}\label{Figure2}
\end{figure}
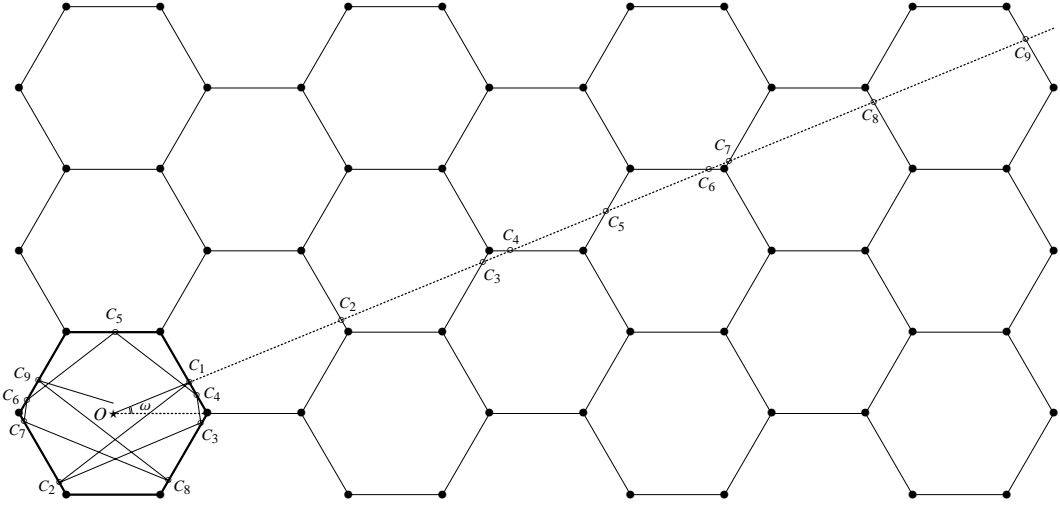

Our original motivation for considering this problem comes from
the study of the exit time of the linear motion with specular
cushion collisions on a hexagonal (open) billiard table with
(small) circular open pockets of radius $\eps$ removed from its
corners (see Figure \ref{Figure2}). The starting remark here is
that, after unfolding the hexagon to a honeycomb in $\R^2$, one
can deform the later to $\Z^2_{(3)}$. This process converts the
problem on the hexagonal billiard with pockets into one concerning
the free path length of a Lorentz gas in $\R^2$ with small
identical ellipses centered at the points from $\Z^2_{(3)}$ as
scatterers (see Figure \ref{Figure8}).

Let $\tau^{\,\operatorname{hex}}_\eps (\omega)$ denote the free
path length in the hexagonal billiard with discs of radius $\eps$
removed from the corners and motion starting at the center, and let
\begin{equation*}
{\mathbb P}^{\,\operatorname{hex}}_\eps (\lambda)=\frac{1}{2\pi}
\bigg| \bigg\{ \omega \in [0,2\pi]:
\tau^{\,\operatorname{hex}}_\eps (\omega)>\frac{\lambda}{\eps}
\bigg\} \bigg|
\end{equation*}
denote the repartition function of $\eps
\tau^{\,\operatorname{hex}}_\eps$. We prove

\begin{theorem}\label{T2}
For every $\delta >0$, as $\eps\rightarrow 0^+$,
\begin{equation*}
{\mathbb P}^{\,\operatorname{hex}}_{\eps} (\lambda) = G_3 \left(
\frac{2\lambda}{\sqrt{3}}\right)+O_\delta \left(
\eps^{\frac{1}{8}-\delta}\right),\qquad \lambda >0.
\end{equation*}
\end{theorem}

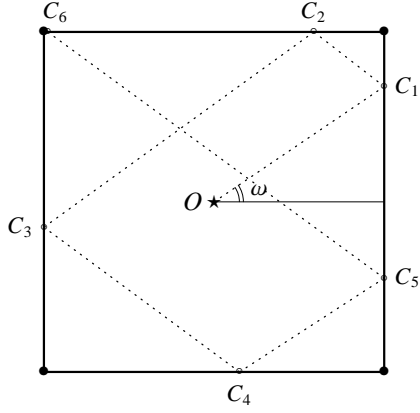
\begin{figure}[ht]
\begin{center}
\unitlength 0.52mm
\begin{picture}(100,95)(0,-3)
\put(87,73){\makebox(0,0){{\tiny $\circ$}}}
\put(69,87){\makebox(0,0){{\tiny $\circ$}}}
\put(0,37){\makebox(0,0){{\tiny $\circ$}}}
\put(50,0){\makebox(0,0){{\tiny $\circ$}}}
\put(87,24){\makebox(0,0){{\tiny $\circ$}}}
\put(1,87){\makebox(0,0){{\tiny $\circ$}}}

\put(0,0){\makebox(0,0){{\small $\bullet$}}}
\put(0,87){\makebox(0,0){{\small $\bullet$}}}
\put(87,87){\makebox(0,0){{\small $\bullet$}}}
\put(87,0){\makebox(0,0){{\small $\bullet$}}}
\put(43.5,43.5){\makebox(0,0){{\small $\star$}}}

\put(55,46.5){\makebox(0,0){\small $\omega$}}
\put(43.5,43.5){\arc{15}{-0.6}{0}}
\put(43.5,43.5){\arc{13}{-0.6}{0}}

\put(93,73){\makebox(0,0){{\small $C_1$}}}
\put(69,92){\makebox(0,0){{\small $C_2$}}}
\put(-6,37){\makebox(0,0){{\small $C_3$}}}
\put(50,-6){\makebox(0,0){{\small $C_4$}}}
\put(93,24){\makebox(0,0){{\small $C_5$}}}
\put(3,92){\makebox(0,0){{\small $C_6$}}}
\put(38,43.5){\makebox(0,0){$O$}}

\thinlines
\path(43.5,43.5)(87,43.5)
\dottedline{2}(43.5,43.5)(87,73)(69,87)(0,37)(50,0)(87,24)(1,87)

\thicklines \path(0,0)(87,0)(87,87)(0,87)(0,0)

\end{picture}
\end{center}
\caption{\it \small The free path in a
square billiard with trajectory starting at the center}\label{Figure3}
\end{figure}

Scaling $\eps$ to $\eps \sqrt{2}$ one can apply Theorem \ref{T1}
(ii) with $\ell=2$ to estimate the repartition function ${\mathbb
P}^\square_\eps (\lambda)=\frac{1}{2\pi} \vert \{ \omega \in
[0,2\pi]: \tau_\eps^\square (\omega)>\frac{\lambda}{\eps}\}\vert$
of the free path length $\tau_\eps^\square (\omega)$ of a billiard
in the unit square with pockets of radius $\eps$ at the corners
and trajectory starting at the center (see Figure \ref{Figure3}),
getting

\begin{theorem}\label{T3}
For every $\delta >0$, as $\eps \rightarrow 0^+$,
\begin{equation*}
{\mathbb P}^\square_\eps (\lambda)=G_2 \left(
\frac{\lambda}{\sqrt{2}}\right)+O_\delta \Big(
\eps^{\frac{1}{8}-\delta}\Big),\qquad \lambda >0.
\end{equation*}
\end{theorem}

Theorem \ref{T3} should be compared with the situation
where the initial position is at one of the four
vertices, and where the limiting distribution has compact support
\cite{BGZ1,BGZ2}. It is not clear
whether these methods would directly extend to other concrete initial positions,
such as $(\frac{1}{n},0)$, $n\in \N$. However, it looks likely that further
refinements could lead to ``space-phase average" results similar to those from
\cite{BZ2}. The main difficulty seems to arise from the increasing complexity
and number of cases that need to be analyzed in detail, leading to
integrals in the main terms of the asymptotic formula which are manifestly more
intricate than in the case of the square.


\section{The contribution of consecutive Farey fractions
$\gamma<\gamma^\prime$ with $\gamma,\gamma^\prime \in
\FF_{Q}\setminus \Fl$}
The integer part of a real number $x$ is denoted by $[x]$.
Let $\eps >0$ and $Q=\big[\frac{1}{\eps}\big]$. 
Denote by $\FF_Q$ the set of Farey fractions
$\gamma=\frac{a}{q}$ in lowest terms with $0<a\leqslant q\leqslant
Q$. For any interval $I\subseteq [0,1]$, set
$\FF_{I,Q}=I\cap\FF_Q$. Denote by $\FF_Q^{(\ell)}$ the set of
Farey fractions $\gamma=\frac{a}{q}\in\FF_Q$ with $\ell \mid
(q-a)$, and set $\FF^{(\ell)}_{I,Q}=I\cap \FF_Q^{(\ell)}$. Set
also $\FF^{(\ell)}=\cup_{Q= 1}^\infty \FF_Q^{(\ell)}$. 
It is well known that if $\gamma=\frac{a}{q}<\gamma^\prime=\frac{a^\prime}{q^\prime}$
are consecutive elements in $\FF_Q$ if and only if
\begin{equation*}
a^\prime q-aq^\prime =1\quad \mbox{\rm and}\quad
q+q^\prime >Q\geqslant \max \{ q,q^\prime \}.
\end{equation*}
In particular we have
\begin{equation}\label{2.1}
\eps (q+q^\prime)\geqslant \eps (Q+1) >1 \geqslant \max \{ \eps
q,\eps q^\prime \} .
\end{equation}

For convenience consider
\begin{equation*}
\tilde{\mathbb G}_{\ell,I,\eps} (\lambda):=\big| \{ \omega \in
\arctan I: q_{\ell,\eps} (\omega)>\lambda Q\} \big| = {\mathbb
G}_{\ell,I,\eps} \bigg( \lambda\eps \bigg[
\frac{1}{\eps}\bigg]\Bigg),\quad \lambda >0.
\end{equation*}
For any interval $I\subseteq [0,1]$ with $\vert I\vert \asymp Q^{-c}$,
$0<c<1$, we will prove a formula of type
\begin{equation*}
\tilde{\mathbb G}_{\ell,I,\eps} (\lambda)=c_I G_\ell
(\lambda)+O_\delta \left( \eps^{-\delta+\theta (c,c_1)}\right),
\end{equation*}
which will immediately imply \eqref{1.1}, because $\eps
\big[ \frac{1}{\eps}\big]=1+O(\eps)$ and for any compact
set $K\subseteq \R_+$ there is $C_K >0$ such that
\begin{equation}\label{2.2}
\vert G_\ell (x^\prime)-G_\ell (x^{\prime \prime})\vert \leqslant
C_K \vert x^\prime -x^{\prime \prime} \vert ,\qquad
x^\prime,x^{\prime \prime}\in K.
\end{equation}

Given $\gamma <\gamma^\prime$ consecutive elements in $\FF_Q$,
denote $t_0=\gamma^\prime -\frac{\eps}{q^\prime}=\frac{a^\prime
-\eps}{q^\prime}$, $u_0=\gamma+\frac{\eps}{q}=\frac{a+\eps}{q}$.
Employing \eqref{2.1} and
\begin{equation}\label{2.3}
t_0 -\gamma =\frac{1-\eps q}{qq^\prime},\qquad u_0 -t_0 =\frac{\eps
(q+q^\prime)-1}{qq^\prime},\qquad \gamma^\prime -u_0 =\frac{1-\eps
q^\prime}{qq^\prime},
\end{equation}
we find $\gamma\leqslant t_0 <u_0 \leqslant \gamma^\prime$. In
particular if $\gamma<\gamma^\prime <\gamma^{\prime \prime}$ are
consecutive elements in $\FF_Q$, then $\gamma +\frac{\eps}{q}
\leqslant \gamma^\prime <\gamma^{\prime \prime}
-\frac{\eps}{q^{\prime \prime}}$. Therefore the intervals $[\gamma
-\frac{\eps}{q},\gamma+\frac{\eps}{q}]$, $\gamma=\frac{a}{q} \in
\FF_Q$, cover the interval $[0,1]$ in such a way that every
element in $[0,1]$ belongs to at most two of these intervals. As a
result any trajectory with slope $\tan\omega \in
(\gamma,\gamma^\prime)$ will intersect, as in the case of the
square lattice \cite{BGZ1,BGZ2}, one of the scatterers
$(q,a)+V_\eps$ or $(q^\prime,a^\prime)+V_\eps$.

Since $a^\prime q-aq^\prime =1$, only two situations can occur here:
$\gamma \notin \FF^{(\ell)}$ and $\gamma^\prime \notin \FF^{(\ell)}$,
respectively $\gamma\in \FF^{(\ell)}$ or $\gamma^\prime \in \FF^{(\ell)}$.
The later will be discussed in Section 3. In the remainder of this section
we assume that $\gamma \notin \FF^{(\ell)}$ and $\gamma^\prime \notin \FF^{(\ell)}$,
situation where the
horizontal free path is given (see Figure \ref{Figure4}) by
\begin{equation*}
q_{\ell,\eps} (\omega)=\begin{cases}
q & \mbox{\rm if $\ \gamma <\tan \omega < t_0,$} \\
\min\{ q,q^\prime \} & \mbox{\rm if $\ t_0 <\tan\omega < u_0,$} \\
q^\prime & \mbox{\rm if $\ u_0 <\tan\omega < \gamma^\prime.$}
\end{cases}
\end{equation*}

\begin{figure}[ht]
\begin{center}
\unitlength 0.6mm
\begin{picture}(160,35)(0,-5)

\path(18,10)(18,30) \path(56,10)(56,30) \path(120,10)(120,30)
\path(17,13)(18,10)(19,13) \path(55,13)(56,10)(57,13)
\path(119,13)(120,10)(121,13)

\thinlines \path(0,1)(0,-1) \path(36,1)(36,-1) \path(75,1)(75,-1)
\path(160,1)(160,-1)

\put(0,-6){\makebox(0,0){{$\gamma$}}}
\put(36,-6){\makebox(0,0){{$t_0$}}}
\put(75,-6){\makebox(0,0){{$u_0$}}}
\put(160,-6){\makebox(0,0){{$\gamma^\prime$}}}
\put(0,20){\makebox(0,0){{$q_{\ell,\eps}(\omega)=q$}}}
\put(85,20){\makebox(0,0){{$q_{\ell,\eps}(\omega)=\min\{
q,q^\prime\}$}}}
\put(140,20){\makebox(0,0){{$q_{\ell,\eps}(\omega)=q^\prime$}}}

\Thicklines \path(0,0)(160,0)
\end{picture}
\end{center}
\caption{The horizontal free path when $\gamma,\gamma^\prime \in\FF_{I,Q}\setminus \Fl$}\label{Figure4}
\end{figure}

In this way the contribution of the interval
$[\gamma,\gamma^\prime]$ to $\tilde{\mathbb
G}_{\ell,I,\eps}(\lambda)$ is given by
\begin{equation*}
\begin{cases}
0 & \mbox{\rm if $\ \max \{ q,q^\prime \} \leqslant \lambda Q,$} \\
\arctan \gamma^\prime -\arctan \gamma & \mbox{\rm if $\ \min \{
q,q^\prime \} >\lambda Q,$} \\
\arctan \gamma^\prime -\arctan u_0 & \mbox{\rm if $\ q
\leqslant \lambda Q < q^\prime,$} \\
\arctan t_0 -\arctan \gamma & \mbox{\rm if $\ q^\prime \leqslant
\lambda Q< q.$}
\end{cases}
\end{equation*}
This contribution is zero whenever $\lambda\geqslant 1$, so we
next assume $0<\lambda <1$.

Using \eqref{2.3}, the estimates
\begin{equation}\label{2.4}
\arctan (x+h)-\arctan x=\frac{h}{1+x^2}+O(h^2),\qquad
\frac{1}{1+\gamma^2}-\frac{1}{1+\gamma^{\prime 2}}\leqslant
2(\gamma^\prime -\gamma)=\frac{2}{qq^\prime},
\end{equation}
and the inequality
\begin{equation*}
\sum_{\gamma\in\FF_Q} \frac{1}{q^2q^{\prime 2}}\leqslant
\sum_{\gamma\in\FF_Q} \frac{1}{Qqq^\prime}=\frac{1}{Q},
\end{equation*}
we infer that the contribution to $\tilde{\mathbb
G}_{\ell,I,\eps}(\lambda)$ of all intervals
$[\gamma,\gamma^\prime]\subseteq I$ with $\gamma <\gamma^\prime$
consecutive elements in $\FF_Q$ and $\gamma,\gamma^\prime \notin
\FF^{(\ell)}$ is given by
\begin{equation*}
{\mathbb G}^{(1)}_{\ell,I,\eps}(\lambda)
=A_{I,Q}(\lambda)+B_{I,Q}(\lambda)+C_{I,Q}(\lambda) +O(\eps),
\end{equation*}
with

\begin{equation}\label{2.5}
A_{I,Q}(\lambda) =\sum\limits_{\substack{\gamma,\gamma^\prime
\in\FF_{I,Q}\setminus \FF^{(\ell)} \\ \min \{ q,q^\prime\}>
\lambda Q}} \frac{1}{qq^\prime}\cdot\frac{1}{1+\gamma^2} ,
\end{equation}
\begin{equation}\label{2.6}
B_{I,Q} (\lambda) = \hspace{-8pt}
\sum\limits_{\substack{\gamma,\gamma^\prime\in\FF_{I,Q}\setminus
\FF^{(\ell)} \\ q\leqslant \lambda Q < q^\prime}} \frac{1-\eps
q^\prime}{qq^\prime}\cdot \frac{1}{1+\gamma^{\,\prime\, 2}} ,\qquad
C_{I,Q} (\lambda)  = \hspace{-8pt}
\sum\limits_{\substack{\gamma,\gamma^\prime\in\FF_{I,Q}\setminus
\FF^{(\ell)} \\ q^\prime \leqslant \lambda Q< q}} \frac{1-\eps
q}{qq^\prime}\cdot \frac{1}{1+\gamma^2}.
\end{equation}

Since there are no consecutive elements in $\FF_Q$ which belong
both to $\Fl$ we have
\begin{equation}\label{2.7}
A_{I,Q}(\lambda)=A^+_{I,Q}(\lambda)-A_{I,Q}^{-,1} (\lambda)
-A_{I,Q}^{-,2} (\lambda),
\end{equation}
with
\begin{equation*}
\begin{split}
A_{I,Q}^+ (\lambda) & =\sum_{\substack{\gamma,\gamma^\prime \in \FF_{I,Q} \\
\min\{ q,q^\prime\}>\lambda Q}} \frac{1}{qq^\prime}\cdot \frac{1}{1+\gamma^2},\\
A_{I,Q}^{-,1} (\lambda) & =\sum_{\substack{\gamma\in \FF_{I,Q}^{(\ell)},\,
\gamma^\prime \in \FF_{I,Q}\\ \min\{ q,q^\prime\} >\lambda Q}}\frac{1}{qq^\prime}\cdot \frac{1}{1+\gamma^2},
\qquad A_{I,Q}^{-,2} (\lambda) =\sum_{\substack{\gamma\in\FF_{I,Q} ,\,
\gamma^\prime\in \FF_{I,Q}^{(\ell)} \\ \min\{ q,q^\prime\} >\lambda Q}} \frac{1}{qq^\prime}\cdot \frac{1}{1+\gamma^2},
\end{split}
\end{equation*}
and respectively
\begin{equation}\label{2.8}
B_{I,Q}(\lambda)=B_{I,Q}^+ (\lambda)-B_{I,Q}^{-,1}(\lambda)-B_{I,Q}^{-,2}(\lambda),\qquad
C_{I,Q}(\lambda)=C_{I,Q}^+ (\lambda)-C_{I,Q}^{-,1}(\lambda)-C_{I,Q}^{-,2}(\lambda),
\end{equation}
with
\begin{equation*}
\begin{split}
& B_{I,Q}^+ (\lambda) = \sum_{\substack{\gamma,\gamma^\prime \in \FF_{I,Q} \\
q\leqslant \lambda Q<q^\prime}}\frac{1-\eps
q^\prime}{qq^\prime}\cdot \frac{1}{1+\gamma^{\,\prime\, 2}},\qquad \quad \ \
C_{I,Q}^+(\lambda) =\sum\limits_{\substack{\gamma,\gamma^\prime\in\FF_{I,Q}
 \\ q^\prime \leqslant \lambda Q< q}} \frac{1-\eps
q}{qq^\prime}\cdot \frac{1}{1+\gamma^2},
\\
& B_{I,Q}^{-,1} (\lambda) = \sum_{\substack{\gamma\in\FF_{I,Q}^{(\ell)},\,\gamma^\prime\in\FF_{I,Q} \\
q\leqslant \lambda Q<q^\prime}}\frac{1-\eps
q^\prime}{qq^\prime}\cdot \frac{1}{1+\gamma^{\,\prime\, 2}},\qquad
C_{I,Q}^{-,1}(\lambda) =\sum\limits_{\substack{\gamma\in\FF_{I,Q}^{(\ell)},\,\gamma^\prime\in\FF_{I,Q}
 \\ q^\prime \leqslant \lambda Q< q}} \frac{1-\eps
q}{qq^\prime}\cdot \frac{1}{1+\gamma^2}, \\ &
B_{I,Q}^{-,2} (\lambda) = \sum_{\substack{\gamma\in\FF_{I,Q},\,\gamma^\prime \in \FF_{I,Q}^{(\ell)} \\
q\leqslant \lambda Q<q^\prime}}\frac{1-\eps
q^\prime}{qq^\prime}\cdot \frac{1}{1+\gamma^{\,\prime\, 2}},\qquad
C_{I,Q}^{-,2} (\lambda) =\sum\limits_{\substack{\gamma,\gamma^\prime\in\FF_{I,Q}
 \\ q^\prime \leqslant \lambda Q< q}} \frac{1-\eps
q}{qq^\prime}\cdot \frac{1}{1+\gamma^2}.
\end{split}
\end{equation*}

\begin{lemma}\label{L2.1}
For any function $V\in C^1[0,N]$ with total variation $T_0^N V$
and $C(\ell)$ as in \eqref{1.2},
\begin{equation*}
\sum\limits_{\substack{1\leqslant k\leqslant N \\ \gcd
(\ell,k)=1}} \frac{\varphi (k)}{k}\, V(k)=C(\ell) \int_0^N V
+O_\ell\Big( (\| V\|_\infty+T_0^N V)\ln N\Big).
\end{equation*}
\end{lemma}

\begin{proof}
The left-hand side $L_{\ell,N}$ can be expressed as
\begin{equation*}
\sum\limits_{\substack{1\leqslant k\leqslant N \\ \gcd
(\ell,k)=1}} \sum\limits_{d\mid k} \frac{\mu (d)}{d}\, V(k).
\end{equation*}
Writing $k=dk^\prime$ with $1\leqslant k^\prime \leqslant
[\frac{N}{d}]$ and $\gcd (\ell,d)=\gcd (\ell,k^\prime)=1$,
denoting $V_n(x)=V(nx)$, $\sigma_0(\ell):=\# \{ d\geqslant 1:
d\mid \ell\}$, and using M\" obius summation
(as in \cite[Lemma 2.2]{BCZ}) we have
\begin{equation*}
\begin{split}
L_{\ell,N} & =\sum\limits_{\substack{1\leqslant d\leqslant N \\
\gcd (\ell,d)=1}} \frac{\mu (d)}{d}
\sum\limits_{\substack{1\leqslant k^\prime \leqslant [N/d]
\\ \gcd (\ell,k^\prime)=1}} V_d (k^\prime) \\ &
=\sum\limits_{\substack{1\leqslant d\leqslant N \\
\gcd (\ell,d)=1}} \frac{\mu (d)}{d} \left( \frac{\varphi
(\ell)}{\ell}\int_0^{[N/d]} V_d
+O\bigg( \Big(\| V_d\|_\infty+T_0^{[N/d]} V_d\Big)\sigma_0 (\ell)\right) \\
& =\frac{\varphi (\ell)}{\ell} \sum\limits_{\substack{1\leqslant d\leqslant N \\
\gcd (\ell,d)=1}} \frac{\mu (d)}{d} \cdot \frac{1}{d}
\left( \int_0^N V+O\big( d\| V\|_\infty \big)\right)+O_\ell \Bigg(
\big( \| V\|_\infty +T_0^N V\big) \sum_{d=1}^N \frac{1}{d} \Bigg)                \\
& =\frac{\varphi (\ell)}{\ell}\sum\limits_{\substack{1\leqslant d\leqslant N \\
\gcd (\ell,d)=1}} \frac{\mu (d)}{d^2}\int_0^N V+O\big( \| V\|_\infty \ln N\big)
+O_\ell \Big( \big( \| V\|_\infty +T_0^N V\big)\ln N \Big)
\\ & =
\frac{\varphi (\ell)}{\ell}\Bigg( \sum\limits_{\substack{d\geqslant 1 \\
\gcd (\ell,d)=1}} \hspace{-5pt} \frac{\mu (d)}{d^2}+O\bigg(
\sum_{d>N} \frac{1}{d^2}\bigg)\Bigg) \int_0^N \hspace{-5pt}
V+O_\ell\Big( \big(\| V\|_\infty+T_0^N V\big)\ln N \Big) \\
& =C(\ell)  \int_0^N V+O_\ell \Big( (\| V\|_\infty+T_0^N V)\ln N
\Big),
\end{split}
\end{equation*}
with
\begin{equation*}
C(\ell)=\frac{\varphi(\ell)}{\ell}
\sum\limits_{\substack{d\geqslant 1 \\ \gcd (\ell,d)=1}} \frac{\mu
(d)}{d^2}=\frac{\varphi(\ell)}{\ell} \prod\limits_{\substack{p\nmid \ell
\\ p\operatorname{prime}}} \left(
1-\frac{1}{p^2}\right)=\frac{\varphi (\ell)}{\zeta (2)\ell}
\prod\limits_{\substack{p\mid \ell \\ p \operatorname{prime}}}
\left( 1-\frac{1}{p^2}\right)^{-1},
\end{equation*}
which gives the desired estimate.
\end{proof}

\begin{lemma}\label{L2.2}
For any function $V\in C^1 [0,N]$ and any $\delta >0$,
\begin{equation*}
\sum_{n=1}^N \frac{\varphi (\ell n)}{n}\ V(n)=\ell C(\ell)
\int_0^N V +O_{\ell,\delta} \left( (\| V\|_\infty +T_0^N V)N^\delta\right).
\end{equation*}
\end{lemma}

\begin{proof}
Let $\ell=p_1^{\alpha_1}\cdots p_r^{\alpha_r}$ with
$p_1,\ldots,p_r$ distinct primes and $\alpha_1,\ldots,\alpha_r \in
\N$ (so $r=\omega (\ell)$, the number of prime divisors of
$\ell$). Writing $n=p_1^{k_1}\cdots p_r^{k_r} m$ with $k_i
\geqslant 0$ and $\gcd (\ell,m)=1$, we obtain
\begin{equation}\label{2.9}
\sum\limits_{1\leqslant n\leqslant N} \frac{\varphi (\ell
n)}{n}\ V(n) =\sum\limits_{\substack{k_1,\ldots,k_r\geqslant 0 \\
p_1^{k_1}\cdots p_r^{k_r}\leqslant N}} \prod_{1\leqslant
i\leqslant r} \frac{\varphi (p_i^{k_i+\alpha_i})}{p_i^{k_i}}
\sum\limits_{\substack{1\leqslant m \leqslant
N/(p_1^{k_1}\cdots p_r^{k_r}) \\ \gcd (\ell,m )=1}}
\frac{\varphi (m)}{m}\ V_{p_1^{k_1}\cdots p_r^{k_r}}
(m).
\end{equation}
According to Lemma \ref{L2.1} the inner sum above can be expressed
as
\begin{equation*}
\begin{split}
C(\ell) \int_0^{N/(p_1^{k_1}\cdots p_r^{k_r})}\hspace{-10pt} &
V(p_1^{k_1} \cdots p_r^{k_r} x)\ dx +O_{\ell}\left( (\| V\|_\infty +T_0^N
V)\ln N\right) \\ & = \frac{C(\ell)}{p_1^{k_1} \cdots
p_r^{k_r}} \int_0^N V+ O_{\ell}\left( (\| V\|_\infty +T_0^N V)\ln N\right).
\end{split}
\end{equation*}
Inserting this into \eqref{2.9} and using the fact that the number
of terms in the first sum in \eqref{2.9} is $\ll (\ln N)^{\omega
(\ell)}$ and $\varphi
(p_i^{k_i+\alpha_i})=p_i^{k_i+\alpha_i}(1-p_i^{-1})$ we infer that
the expression in \eqref{2.9} is given by
\begin{equation}\label{2.10}
\ell C(\ell)\prod\limits_{\substack{p\mid \ell \\ p\operatorname{prime}}}
\left( 1-\frac{1}{p}\right)\sum\limits_{\substack{k_1,\ldots,k_r\geqslant 0 \\
p_1^{k_1}\cdots p_r^{k_r}\leqslant N}} \frac{1}{p_1^{k_1}\cdots
p_r^{k_r}} \int_0^N V +O_{\ell} \left( (\| V\|_\infty+T_0^N V)(\ln
N)^{1+\omega(\ell)}\right).
\end{equation}
The statement now follows from \eqref{2.10}, using
\begin{equation*}
\sum_{k_1,\ldots,k_r \geqslant 0} \frac{1}{p_1^{k_1}\cdots
p_r^{k_r}} =\prod\limits_{\substack{p\mid \ell \\
p\operatorname{prime}}} \left( 1-\frac{1}{p}\right)^{-1}
\end{equation*}
and the bound

\begin{equation*}
\sum_{\substack{k_1,\ldots,k_r \geqslant 0 \\ p_1^{k_1}\cdots
p_r^{k_r}>N}} \frac{1}{p_1^{k_1}\cdots p_r^{k_r}}
\ll_{p_1,\ldots,p_r} \frac{(\ln N)^{r-1}}{N}.
\end{equation*}
\end{proof}

The following estimate \cite[Proposition A4]{BZ1} will be employed several times.

\begin{lemma}\label{L2.3}
Assume that $q\geqslant 1$ and $h$ are two given integers, $\II$ and $\JJ$ are intervals
of length less than $q$, and $f:\II \times \JJ \rightarrow \R$ is a $C^1$ function.
Then for any integer $T>1$ and any $\delta >0$
\begin{equation*}
\sum\limits_{\substack{a\in\II,b\in\JJ \\ ab=h\hspace{-5pt}\pmod{q} \\ \gcd (b,q)=1}}
f(a,b)=\frac{\varphi (q)}{q^2} \iint_{\II\times \JJ} f(x,y)\, dx dy +\EE,
\end{equation*}
with
\begin{equation*}
\EE \ll_\delta T^2 \| f\|_\infty q^{\frac{1}{2}+\delta} \gcd (h,q)^{\frac{1}{2}}
+T \| \nabla f\|_\infty q^{\frac{3}{2}+\delta} \gcd (h,q)^{\frac{1}{2}}
+\frac{\| \nabla f\|_\infty\vert \II \vert \vert \JJ\vert}{T},
\end{equation*}
where $\| f\|_\infty$ and $\| \nabla f \|_\infty$ denote the sup-norm
of $f$ and respectively $\big| \frac{\partial f}{\partial x}\big|
+\big| \frac{\partial f}{\partial y}\big|$ on $\II \times\JJ$.
\end{lemma}

Lemma \ref{L2.3} will be typically applied to the following situations:
Let $I$ be a subinterval of $[0,1]$.
For every $q\in [1,Q]$ consider the intervals $\II=qI$ and
$\JJ=J_{\lambda,q}=(\max\{ \lambda Q, Q-q\},Q]$,
and the functions $f_q$ and $g_q$ defined on $qI \times J_{\lambda,q}$ by
\begin{equation*}
f_{q}(u,v):=\frac{1}{qv}\cdot \frac{1}{1+\left(\frac{u}{q}\right)^2},\qquad
g_q(u,v)=\frac{1-\eps v}{qv}\cdot \frac{1}{1+\left(\frac{u}{q}\right)^2}.
\end{equation*}
We clearly have
\begin{equation}\label{2.11}
\| f_{q} \|_\infty \ll_{\lambda} \frac{1}{Qq},\qquad
\| \nabla f_{q} \|_\infty \ll_{\lambda} \frac{1}{Qq^2},\qquad
\| g_q\|_\infty \ll_\lambda \frac{1}{Qq},\qquad \| \nabla g_q\|_\infty \ll_\lambda
\frac{1}{Qq^2}.
\end{equation}

\begin{proposition}\label{P2.1}
For every $\lambda \in (0,1]$ and $c_1\in (0,1)$,
\begin{equation*}
A_{I,Q} (\lambda) =c_I A(\ell)I_1 (\lambda)+O_{\delta,\lambda,\ell} (Q^{\delta+\theta_1(c_1)}),
\end{equation*}
where
\begin{equation*}
I_1(\lambda)=\int_\lambda^1 \frac{1}{x}\, \ln\frac{1}{\max\{ \lambda,1-x\}}\ dx=
\begin{cases} \ln (1-\lambda)\ln\lambda +\int_\lambda^{1-\lambda} \frac{1}{x}\, \ln\frac{1}{1-x}
\ dx & \mbox{if $\lambda\in \big(0,\frac{1}{2}\big],$} \\
(\ln\lambda)^2 & \mbox{if $\lambda \in \big[\frac{1}{2}, 1\big],$}
\end{cases}
\end{equation*}
and $\theta_1 (c_1)=\max\left\{ 2c_1 -\frac{1}{2},-c_1,-1\right\}$.
\end{proposition}

\begin{proof}
There is at most one $\gamma\in\FF_{I,Q}$ with $\gamma^\prime \notin I$.
Since $\frac{1}{\vert I\vert}\ll Q^c$ and $\frac{1}{qq^\prime} \leqslant \frac{1}{Q}$,
the total contribution to the final asymptotic results from
Theorem \ref{T1} resulting from replacing the two conditions
$\gamma,\gamma^\prime \in \FF_{I,Q}$ by $a\in qI$ or by $a^\prime \in q^\prime I$
will be $\ll Q^{-1}$ and respectively $\ll Q^{\frac{1}{8}-1}$, thus negligible.
As a result we shall tacitly do this in formulas \eqref{2.12},
\eqref{2.15}, \eqref{2.19}, \eqref{2.24} and \eqref{2.27}.

Furthermore, the summation constraints in $\gamma=\frac{a}{q}\in \FF_{I,Q}$ and
$\gamma^\prime=\frac{a^\prime}{q^\prime}\in \FF_Q$ will translate,
using standard properties of Farey fractions, into the following constraints on the triplet
$(q,a,q^\prime)$:
\begin{equation*}
\begin{cases}
q\in (\lambda Q,Q],\\
a\in qI,\ \ q^\prime \in J_{\lambda,q},\ \
\gcd (q^\prime ,q)=1,\ \ aq^\prime =-1\hspace{-5pt} \pmod{q}.
\end{cases}
\end{equation*}
Summing first over $q$ and then after integer pairs $(a,q^\prime)\in qI\times J_{\lambda,q}$
as above (note that the number of such pairs is $\leqslant q\vert I\vert$) we infer
\begin{equation}\label{2.12}
A_{I,q}^+(\lambda) =\sum_{q\in (\lambda Q,Q]}
\sum_{\substack{(a,q^\prime)\in qI\times J_{\lambda,q} \\ \gcd (q^\prime,q)=1 \\
aq^\prime \equiv -1 \hspace{-5pt}\pmod{q}}} f_q (a,q^\prime)+O\left( \frac{1}{Q}\right).
\end{equation}
Applying Lemma \ref{L2.3} with $T=[Q^{c_1}]$ and \eqref{2.11}, the inner sum in \eqref{2.12}
can be expressed as
\begin{equation}\label{2.13}
\begin{split}
\frac{\varphi (q)}{q^2} & \int_{J_{\lambda,q}} \frac{dv}{qv} \int_{qI}
\frac{du}{1+\left( \frac{u}{q}\right)^2} \ +\
O_{\delta,\lambda}\left( Q^{2c_1} q^{\frac{1}{2}+\delta} Q^{-2}+
Q^{c_1} q^{\frac{3}{2}+\delta} Q^{-3} +q^2 Q^{-3}Q^{-c_1}\right) \\
& =c_I \frac{\varphi (q)}{q^2} \, \ln
\frac{1}{\max\big\{ \lambda ,1-\frac{q}{Q}\big\}}+
O_{\delta,\lambda} \left( Q^{2c_1-\frac{3}{2}+\delta}+Q^{-1-c_1}\right).
\end{split}
\end{equation}
Summing over $q\in (\lambda Q,Q]$ in \eqref{2.13} and employing \eqref{2.12}
and \cite[Lemma 2.3]{BCZ} we find
\begin{equation}\label{2.14}
A_{I,Q}^+ (\lambda) = \frac{c_I}{\zeta(2)} \, I_1 (\lambda)+O_{\delta,\lambda}
\left(Q^{\delta+\theta_1(c_1)}\right).
\end{equation}

To estimate $A_{I,Q}^{-,1}(\lambda)$, note first that $\gcd
(\ell,q)=1$ because $\ell \mid q-a$ and $\gcd (q,a)=1$. As a
result, putting $w=q-a=\ell u$, $v=q^\prime$ and using $a^\prime
q-aq^\prime=1$ and $q^\prime >Q-q$, we find (with $\bar{\ell}$
denoting the multiplicative inverse of $\ell\pmod{q}$)
\begin{equation}\label{2.15}
\begin{split}
A_{I,Q}^{-,1}(\lambda)  &
=\sum_{\substack{\lambda Q <q\leqslant Q \\ \gcd (\ell,q)=1}}
\sum_{\substack{(w,v)\in q(1-I)\times J_{\lambda,q} \\ \ell \mid w,\, \gcd (v,q)=1 \\
wv\equiv 1 \hspace{-5pt} \pmod{q}}}
\frac{1}{qv}\cdot \frac{1}{1+\left( \frac{q-w}{q}\right)^2} +O\left( \frac{1}{Q}\right)  \\ &
=\sum_{\substack{\lambda Q <q\leqslant Q \\ \gcd (\ell,q)=1}}
\sum_{\substack{(u,v)\in (q/\ell)(1-I)\times \in J_{\lambda,q} \\ \gcd (v,q)=1 \\
uv=\bar{\ell} \hspace{-5pt}\pmod{q}}}
f_{q} (q-\ell u,v)+O\left( \frac{1}{Q}\right).
\end{split}
\end{equation}
By Lemma \ref{L2.3} (with $T=[Q^{c_1}]$) and
\eqref{2.11} the inner sum in \eqref{2.15} can be expressed as
\begin{equation}\label{2.16}
\begin{split}
\frac{\varphi (q)}{q^2} & \int_{J_{\lambda,q}} \frac{dv}{qv} \int_{\frac{q}{\ell} (1-I)}
\frac{du}{1+\left( \frac{q-\ell u}{q}\right)^2} \ +\
O_{\delta,\lambda,\ell} \left( Q^{2c_1-\frac{3}{2}+\delta}+Q^{-1-c_1}\right) \\
& =c_I \frac{\varphi (q)}{\ell q^2} \, \ln \frac{1}{\max\big\{ \lambda ,1-\frac{q}{Q}\big\}}+
O_{\delta,\lambda,\ell} \left( Q^{2c_1-\frac{3}{2}+\delta}+Q^{-1-c_1}\right).
\end{split}
\end{equation}
Summing over $q$ in \eqref{2.16} and employing \eqref{2.15} and
Lemma \ref{L2.1} we find
\begin{equation}\label{2.17}
A_{I,Q}^{-,1} (\lambda) \ =\frac{c_IC(\ell)}{\ell} \, I_1 (\lambda)+O_{\delta,\lambda,\ell}
\left(Q^{\delta+\theta_1(c_1)}\right).
\end{equation}

To estimate $A_{I,Q}^{-,2}$ we use the second inequality in \eqref{2.4} and
$\sum_{\gamma\in \FF_{I,Q}} \frac{1}{qq^\prime} \leqslant 1$ to infer
\begin{equation}\label{2.18}
A_{I,Q}^{-,2}=\sum_{\substack{\gamma\in \FF_{I,Q},\, \gamma^\prime \in \FF_{I,Q}^{(\ell)} \\
\min\{ q,q^\prime \} >\lambda Q}} \frac{1}{qq^\prime}\cdot \frac{1}{1+\gamma^{\,\prime\, 2}}\ +\
O\left( \frac{1}{Q}\right).
\end{equation}
We now proceed as for $A_{I,Q}^{-,1}$, noting that $\gcd (\ell,q^\prime)=1$ because $\ell
\mid q^\prime -a^\prime$ and $\gcd (q^\prime,a^\prime)=1$, and getting as above
(with $w=q^\prime -a^\prime=\ell u$, $v=q$)
\begin{equation}\label{2.19}
\begin{split}
A_{I,Q}^{-,2}(\lambda) & =\sum_{\substack{q^\prime \in (\lambda Q,Q] \\
\gcd (\ell,q^\prime)=1}} \sum_{\substack{(w,v)\in q^\prime (1-I)\times J_{\lambda,q^\prime} \\
\ell \vert w,\, \gcd (v,q^\prime)=1
\\ wv=-1\hspace{-5pt}\pmod{q^\prime}}} \frac{1}{q^\prime v}\cdot
\frac{1}{1+\left( \frac{q^\prime -w}{q^\prime}\right)^2}+O\left(\frac{1}{Q}\right) \\ &
=\sum_{\substack{q^\prime \in (\lambda Q,Q] \\ \gcd
(\ell,q^\prime)=1}} \sum_{\substack{(u,v)\in (q^\prime/\ell)(1-I)\times J_{\lambda,q^\prime} \\
\gcd (v,q^\prime)=1 \\ uv=-\bar{\ell}\hspace{-5pt}
\pmod{q^\prime}}} f_{q^\prime}(q^\prime -\ell u,v)+O\left(\frac{1}{Q}\right)
\\ & =\sum_{\substack{q^\prime \in (\lambda Q,Q] \\ \gcd (\ell,q^\prime)=1}} \left(
\frac{\varphi (q^\prime)}{q^{\prime\,2}} \int_{J_{\lambda,q^\prime}} \frac{dv}{q^\prime v}
\int_{\frac{q^\prime}{\ell} (1-I)} \frac{du}{1+\big(\frac{q^\prime -\ell u}{q^\prime}\big)^2} \ +
\ O_{\delta,\lambda,\ell} \left( Q^{2c_1-\frac{3}{2}+\delta}+Q^{-1-c_1}\right) \right)\\ &
=\frac{c_I}{\ell} \sum_{\substack{q^\prime \in (\lambda Q,Q] \\ \gcd (\ell,q^\prime)=1}}
\frac{\varphi (q^\prime)}{q^{\prime\,2}}\, \ln
\frac{1}{\max\big\{ \lambda,1-\frac{q^\prime}{Q}\big\}}\
+\ O_{\delta,\ell} \left( Q^{2c_1-\frac{1}{2}+\delta}+Q^{-c_1}\right) \\
& =\frac{c_I C(\ell)}{\ell}\ I_1 (\lambda)+O_{\delta,\lambda,\ell}
\left( Q^{\delta+\theta_1 (c_1)}\right).
\end{split}
\end{equation}

The estimate for $A_{I,Q}(\lambda)$ follows from \eqref{2.7}, \eqref{2.14}, \eqref{2.17}
and \eqref{2.18}.
\end{proof}

\begin{proposition}\label{P2.2}
For every $\lambda \in (0,1]$ and $c_1\in (0,1)$,
\begin{equation*}
B_{I,Q} (\lambda) =c_I A(\ell)I_2 (\lambda)+O_{\delta,\lambda,\ell}
\left(Q^{\delta+\theta_1(c_1)}\right)\quad \mbox{and}
\quad C_{I,Q} (\lambda) =c_I A(\ell)I_2 (\lambda)+O_{\delta,\lambda,\ell}
\left(Q^{\delta+\theta_1(c_1)}\right),
\end{equation*}
where
\begin{equation*}
I_2(\lambda) =\int_{\max\{ \lambda,1-\lambda\}}^1 \frac{1-x}{x}\, \ln\frac{\lambda}{1-x}\ dx.
\end{equation*}
\end{proposition}

\begin{proof}
Using the argument leading to \eqref{2.18} we see that
\begin{equation*}
B_{I,Q}^+(\lambda)=\tilde{B}_{I,Q}^+(\lambda)+O\bigg(\frac{1}{Q}\bigg)
=\sum_{\substack{\gamma,\gamma^\prime\in\FF_{I,Q} \\ q\leqslant
\lambda Q<q^\prime}} \frac{1-\eps q^\prime}{qq^\prime}\cdot
\frac{1}{1+\gamma^2} +O\left( \frac{1}{Q}\right).
\end{equation*}
Using customary properties of Farey fractions we infer (setting $u=q-a$, $v=q^\prime$)
\begin{equation*}
\tilde{B}^+_{I,Q}(\lambda)=\sum_{q \leqslant \lambda Q}
\sum_{\substack{(u,v)\in q(1-I)\times J_{\lambda,q}\\
uv\equiv 1 \hspace{-5pt} \pmod{q}}} g_{q}(q-u,v)+O\left( \frac{1}{Q}\right).
\end{equation*}
Applying Lemma \ref{L2.3} to $g_{q}$ and $T=[Q^{c_1}]$ to the inner sum above and
using \eqref{2.11} we find
\begin{equation*}
\begin{split}
B_{I,Q}^+(\lambda) & =\sum_{q \leqslant \lambda Q} \left(
\frac{\varphi (q)}{q^{2}} \int_{q(1-I)}
\frac{du}{1+\frac{(q-u)^2}{q^{2}}} \int_{J_{\lambda,q}}
\frac{1-\eps v}{qv}\ dv + O_{\delta}\Bigg(
\frac{Q^{2c_1}q^{\frac{1}{2}+\delta}}{Qq}+
\frac{Q^{c_1}q^{\frac{3}{2}+\delta}}{Qq^{2}} +
\frac{q^{2} \frac{1}{Qq^{2}}}{Q^{c_1}}\Bigg)\right)  \\
& =c_I \sum_{q\leqslant \lambda Q} \frac{\varphi (q)}{q} \, V(q)+
O_\delta \left( Q^{\max \{ 2c_1-\frac{1}{2}+\delta,-c_1\}}\right),
\end{split}
\end{equation*}
where
\begin{equation*}
V(n)=\frac{1}{n} \int_{\max\{ \lambda,1-\frac{n}{Q}\}}^1 \frac{1-y}{y}\ dy,
\qquad n\in [1,\lambda Q].
\end{equation*}
Using

\begin{equation}\label{2.20}
\| V\|_\infty \leqslant \frac{1}{\lambda Q} \qquad \mbox{\rm and} \qquad
T_0^{\lambda Q} V\ll_\lambda \frac{1}{Q},
\end{equation}
and applying M\" obius summation to $V$ (e.g. \cite[Lemma 2.3]{BCZ}) and Tonelli's
theorem we find
\begin{equation}\label{2.21}
\begin{split}
B^+_{I,Q}(\lambda) & =\frac{c_I}{\zeta(2)} \int_0^{\lambda Q} \frac{du}{u}
\int_{\max\{ \lambda,1-\frac{u}{Q}\}}^1 \frac{1-y}{y}\ dy
+O_\delta \left( Q^{\delta+\theta_1 (c_1)}\right) \\
& =\frac{c_I}{\zeta(2)} \int_0^\lambda \frac{dx}{x} \int_{\max\{ \lambda,1-x\}}^1
\frac{1-y}{y}\ dy +O_\delta \left( Q^{\delta+\theta_1(c_1)}\right) \\
& =\frac{c_I}{\zeta(2)}\, I_2 (\lambda)+O_\delta \left( Q^{\delta+\theta_1 (c_1)}\right).
\end{split}
\end{equation}

The condition $\ell \mid q-a$ gives $\gcd (\ell,q)=1$. Taking $q-a=\ell u$ and $v=q^\prime$ and
proceeding as in the case of $A_{I,Q}^{-,1}$ from the proof of Proposition \ref{P2.1} we have
\begin{equation}\label{2.22}
\begin{split}
B_{I,Q}^{-,1}(\lambda)  &
=\sum_{\substack{\gamma\in\FF_{I,Q}^{(\ell)},\,
\gamma^\prime\in\FF_{I,Q} \\ q\leqslant \lambda Q<q^\prime}}
\frac{1-\eps q^\prime}{qq^\prime}\cdot\frac{1}{1+\gamma^2} \\ &
=\sum_{\substack{q\leqslant \lambda Q \\ \gcd (\ell,q)=1}}
\sum_{\substack{(u,v)\in (q/\ell)(1-I)\times  J_{\lambda,q}\\
\gcd (v,q)=1 \\ uv \equiv \bar{\ell}\hspace{-5pt}\pmod{q}}}
\frac{1-\eps v}{qv}\cdot \frac{1}{1+\left(\frac{q-\ell u}{q}\right)^2} +O\left(\frac{1}{Q}\right) \\
& =\frac{c_I}{\ell} \sum_{\substack{q\leqslant \lambda Q \\ \gcd
(\ell,q)=1}} \frac{\varphi (q)}{q}\, V(q)+ O_{\delta,\ell}
\left( Q^{\delta+\theta_1 (c_1)}\right).
\end{split}
\end{equation}
Applying Lemma \ref{L2.1} to the last sum in \eqref{2.22} we find
\begin{equation}\label{2.23}
\begin{split}
B_{I,Q}^{-,1} & =\frac{c_I C(\ell)}{\ell} \int_0^{\lambda Q}
\frac{du}{u} \int_{\max\{ \lambda,1-\frac{u}{Q}\}}^1
\frac{1-y}{y}\  dy+O_{\delta,\ell}
\left(Q^{\delta+\theta_1(c_1)}\right) \\ &
=\frac{c_IC(\ell)}{\ell}\ I_2(\lambda)+O_{\delta,\ell}
\left(Q^{\delta+\theta_1(c_1)}\right).
\end{split}
\end{equation}

To estimate $B_{I,Q}^{-,2}$ we first fix $q\in (\lambda Q,Q]$,
then set $u=q-a$, $v=q^\prime$. Since $\ell \mid
q^\prime-a^\prime$ and $\gcd(q^\prime,a^\prime)=1$ we have $\gcd
(\ell,q^\prime)=1$. Moreover, $q^\prime-a^\prime=\frac{q^\prime
u-1}{q}$ is divisible by $\ell$, so $uv\equiv
1\pmod{(\ell q)}$ and we have
\begin{equation}\label{2.24}
B_{I,Q}^{-,2} (\lambda)
=\sum_{\substack{\gamma\in\FF_{I,Q},\,\gamma^\prime\in\FF_{I,Q}^{(\ell)}\\
q\leqslant \lambda Q<q^\prime}} \frac{1-\eps
q^\prime}{qq^\prime}\cdot \frac{1}{1+\left(\frac{a}{q}\right)^2}
=\sum_{q\leqslant \lambda Q} \sum_{\substack{(u,v)\in
q(1-I)\times J_{\lambda,q} \\ \gcd (\ell q,v)=1\\ uv\equiv
1\hspace{-5pt} \pmod{(\ell q)}}} g_q (q-u,v)+O\left( \frac{1}{Q}\right).
\end{equation}
By \eqref{2.24}, Lemma \ref{L2.3} and \eqref{2.11} we infer
\begin{equation}\label{2.25}
\begin{split}
B_{I,Q}^{-,2}(\lambda) & =\sum_{q\leqslant \lambda Q}
\Bigg(\frac{\varphi (\ell q)}{\ell^2 q^2} \ qc_I\,
\frac{1}{q}\int_{J_{\lambda,q}} \frac{1-\frac{v}{Q}}{v}\ dv+
O_\delta \left( Q^{2c_1-1} q^{-\frac{1}{2}+\delta}+Q^{-c_1-1}\right)\Bigg) \\
& =\frac{c_I}{\ell^2} \sum_{q\leqslant \lambda Q}\frac{\varphi (\ell q)}{q}\ V(q)
+O_\delta \left(Q^{\max\{ 2c_1-\frac{1}{2}+\delta,-c_1\}}\right).
\end{split}
\end{equation}
From \eqref{2.25}, Lemma \ref{L2.2} and \eqref{2.20} we infer
\begin{equation}\label{2.26}
B_{I,Q}^{-,2}(\lambda) =\frac{c_IC(\ell)}{\ell} \int_0^{\lambda Q} V+O_{\delta,\ell} \left(
Q^{\delta+\theta_1(c_1)}\right)=\frac{c_I C(\ell)}{\ell}\ I_2(\lambda)+O_{\delta,\ell}
\left( Q^{\delta+\theta_1(c_1)}\right).
\end{equation}
The desired estimate on $B_{I,Q}(\lambda)$ follows from \eqref{2.8}, \eqref{2.21},
\eqref{2.23} and \eqref{2.25}.

In similar fashion one gets
\begin{equation}\label{2.27}
\begin{split}
C_{I,Q}^+(\lambda) & =\sum_{\substack{\gamma,\gamma^\prime\in\FF_{I,Q} \\
q^\prime \leqslant \lambda Q <q}} \frac{1-\eps q}{qq^\prime} \cdot
\frac{1}{1+\gamma^{\,\prime\,2}}
=\sum_{q^\prime\leqslant \lambda Q}\
\sum_{\substack{(u,v)=(q^\prime -a^\prime,q) \in q^\prime (1-I)
\times J_{\lambda,q^\prime}\\ \gcd (v,q^\prime)=1,\, uv\equiv -1
\hspace{-5pt} \pmod{q^\prime}}} \hspace{-10pt} g_{q^\prime} (q^\prime -u,v)
+O\left(\frac{1}{Q}\right) \\ & =\frac{c_I}{\zeta (2)}\ I_2 (\lambda)+O_\delta
\left( Q^{\delta+\theta_1(c_1)}\right),\\ C_{I,Q}^{-,1} (\lambda)
& =\sum_{\substack{\gamma\in\FF_{I,Q}^{(\ell)},\,
\gamma^\prime\in\FF_{I,Q} \\ q^\prime \leqslant \lambda Q <q}}
\frac{1-\eps q}{qq^\prime}\cdot \frac{1}{1+\gamma^{\,\prime\,2}}
 =\sum_{q^\prime\leqslant\lambda Q}\
\sum_{\substack{(u,v)=(q^\prime-a^\prime,q)\in q^\prime (1-I) \times
 J_{\lambda,q^\prime}\\ \gcd (v,q^\prime)=1 \,
uv\equiv -1 \hspace{-5pt}\pmod{(\ell q^\prime)}}} \hspace{-10pt} g_{q^\prime}
(q^\prime -u,v)+O\left(\frac{1}{Q}\right) \\ & =\sum_{q^\prime \leqslant \lambda Q}
\Bigg( \frac{\varphi (\ell q^\prime)}{\ell^2 q^{\prime 2}}\
q^\prime c_I V(q^\prime) +O_\delta \left(
Q^{\max\{ 2c_1-\frac{3}{2}+\delta,-c_1-1\}}\right)\Bigg) \\
& =\frac{c_I C(\ell)}{\ell}\ I_2 (\lambda)+O_{\delta,\lambda,\ell}
\left(Q^{\delta+\theta_1(c_1)}\right), \\ C_{I,Q}^{-,2}(\lambda) &
=\sum_{\substack{\gamma\in\FF_{I,Q},\,\gamma^\prime\in\FF_{I,Q}^{(\ell)} \\
q^\prime \leqslant \lambda Q<q}} \frac{1-\eps q}{qq^\prime} \cdot
\frac{1}{1+\gamma^{\,\prime\,2}}
=\sum_{\substack{q^\prime\leqslant \lambda Q \\ \gcd (\ell,q^\prime)=1}}
\sum_{\substack{(u,v)=\big(\frac{q^\prime -a^\prime}{\ell},q\big)\in\frac{q^\prime}{\ell}(1-I)
\times J_{\lambda,q^\prime}\\ \gcd (v,q^\prime)=1,\, uv\equiv -
\bar{\ell}\hspace{-5pt} \pmod{q^\prime}}}
\hspace{-10pt} g_{q^\prime} (q^\prime -\ell u,v)+O\left(\frac{1}{Q}\right) \\
& =\sum_{\substack{q^\prime\leqslant\lambda Q \\ \gcd (\ell,q^\prime)=1}} \Bigg(
\frac{\varphi (q^\prime)}{q^{\prime 2}} \cdot \frac{q^\prime c_I}{\ell}
\int_{J_{\lambda,q^\prime}} \frac{1-\eps v}{q^\prime v}\ dv +O_{\delta,\ell} \left(
Q^{\max\{ 2c_1-\frac{3}{2}+\delta,-c_1-1\}}\right)\Bigg) \\
& =\frac{c_I C(\ell)}{\ell}\ I_2(\lambda)
+O_{\delta,\lambda,\ell} \left(Q^{\delta+\theta_1(c_1)}\right).
\end{split}
\end{equation}
The desired estimate on $C_{I,Q}(\lambda)$ follows from \eqref{2.8} and \eqref{2.27}.
\end{proof}

\begin{corollary}\label{C2.1}
For every $\lambda >0$ and $\delta >0$,
\begin{equation*}
{\mathbb G}^{(1)}_{\ell,I,\eps}(\lambda)=c_I A(\ell)
G^{(1)}(\lambda)+O_{\delta,\ell} \left( \eps^{-\delta+\theta
(c,c_1)}\right),
\end{equation*}
where
\begin{equation*}
G^{(1)} (\lambda)=\begin{cases} \ln (1-\lambda)\ln
\lambda+\int_\lambda^{1-\lambda} \frac{1}{u}\, \ln \frac{1}{1-u}\
du +2\int_{1-\lambda}^1 \frac{1-u}{u}\, \ln \frac{\lambda}{1-u} \
du & \mbox{if $\ \lambda\in \big(0,
\frac{1}{2}\big],$} \\
(\ln\lambda)^2 +2\int_{\lambda}^1 \frac{1-u}{u}\
\ln\frac{\lambda}{1-u}\ du &
\mbox{if $\ \lambda\in\big[ \frac{1}{2}, 1\big],$} \\
0 & \mbox{if $\ \lambda \in [1,\infty).$}
\end{cases}
\end{equation*}
\end{corollary}

\section{The contribution of consecutive Farey fractions
$\gamma<\gamma^\prime$ with $\gamma\in\FF_{Q}^{(\ell)}$ or
$\gamma^\prime\in\FF_{Q}^{(\ell)}$} Suppose first that
$\gamma<\gamma^\prime$ are consecutive in $\FF_{I,Q}$ and
$\gamma\in \FF_{I,Q}^{(\ell)}$, where again
$Q=\big[\frac{1}{\eps}\big]$. Consider
\begin{equation*}
a_k=ka+a^\prime,\quad q_k=kq+q^\prime,\quad
\gamma_k=\frac{a_k}{q_k},\quad
t_k=\gamma_k-\frac{\eps}{q_k}=\frac{a_k-\eps}{q_k}, \quad
k\geqslant 0,\quad t_{-1}=\gamma^\prime.
\end{equation*}

Inequalities \eqref{2.1} show that
\begin{equation}\label{3.1}
\gamma \ \stackrel{k}{\longleftarrow}\ \gamma_{k+1} <\gamma_k <\cdots
<\gamma_1<\gamma_0=\gamma^\prime \quad \mbox{\rm and} \quad
t_{k+1}<t_k <\gamma_{k+1}<\gamma_k .
\end{equation}

Since $q\equiv a\hspace{-3pt}\pmod{\ell}$ and $a^\prime q-aq^\prime
=1$, we cannot have $q^\prime \equiv a^\prime \hspace{-3pt}
\pmod{\ell}$, and so $\gamma^\prime \in\FF_{I,Q}\setminus\FF^{(\ell)}$.
Moreover, since $\ell \mid (q-a)$ we must have $\gcd (\ell,q)=1$
and $q_k \nequiv a_k \hspace{-3pt} \pmod{\ell}$ for all $k\geqslant
0$.

\begin{figure}[ht]
\begin{center}
\unitlength 0.6mm
\begin{picture}(160,35)(0,-5)

\path(0,0)(17,0) \path(33,0)(160,0) \path(25,10)(25,30)
\path(76,10)(76,30) \path(120,10)(120,30)
\path(24,13)(25,10)(26,13) \path(75,13)(76,10)(77,13)
\path(119,13)(120,10)(121,13)

\thinlines \path(0,1)(0,-1) \path(62,1)(62,-1) \path(90,1)(90,-1)
\path(160,1)(160,-1) \path(17,1)(17,-1) \path(33,1)(33,-1)
\path(48,-1)(48,1)

\put(0,-6){\makebox(0,0){{$\gamma$}}}
\put(90,-6){\makebox(0,0){{$t_0$}}}
\put(62,-6){\makebox(0,0){{$t_1$}}}
\put(17,-6){\makebox(0,0){{$t_{k+1}$}}}
\put(33,-6){\makebox(0,0){{$t_k$}}}
\put(48,-6){\makebox(0,0){{$\gamma_{k+1}$}}}
\put(160,-6){\makebox(0,0){{$\gamma_0=\gamma^\prime$}}}
\put(48,20){\makebox(0,0){{$q_{\ell,\eps}(\omega)=q_{k+1}$}}}
\put(95,20){\makebox(0,0){{$q_{\ell,\eps}(\omega)=q_1$}}}
\put(146,20){\makebox(0,0){{$q_{\ell,\eps}(\omega)=q_0=q^\prime$}}}

\Thicklines \path(17,0)(33,0)
\end{picture}
\end{center}
\caption{The horizontal free path when $\gamma\in \FF_{I,Q}^{(\ell)}$
and $\gamma^\prime \in \FF_{I,Q}\setminus\FF^{(\ell)}$}\label{Figure5}
\end{figure}

The number $K=\big[ \frac{\lambda Q-q^\prime}{q}\big]$ is the
unique integer $K\geqslant 1$ for which $q_K =Kq+q^\prime \leqslant \lambda
Q <q_{K+1}=(K+1)q+q^\prime$. The presence of the sink at $\gamma$,
\eqref{3.1} and $\gamma_k \notin \FF^{(\ell)}$ show that (see also
Figure \ref{Figure5})
\begin{equation*}
q_{\ell,\eps}(\omega)=q_{k+1} \quad \mbox{\rm if
$t_{k+1}<\tan\omega <t_k$, $\ k\geqslant -1$,}
\end{equation*}
and the contribution to ${\mathbb G}_{I,Q}(\lambda)$ of the
interval $[\gamma,\gamma^\prime]$ is given by
\begin{equation}\label{3.2}
\arctan t_K-\arctan \gamma=\frac{1-\eps
q}{q(Kq+q^\prime)}\cdot\frac{1}{1+\gamma^2} +O\left( \frac{1}{q^2
(Kq+q^\prime)^2}\right).
\end{equation}
Since $qq^\prime \geqslant Q$ and $\sum_{\gamma\in\FF_{I,Q}}
\frac{1}{qq^\prime}\leqslant 1$, it follows from \eqref{3.2} that the
total contribution to $\tilde{\mathbb G}_{\ell,I,\eps}(\lambda)$
of intervals $[\gamma,\gamma^\prime]\subseteq I$ with
$\gamma<\gamma^\prime$ consecutive in $\FF_Q$ and
$\gamma\in\FF^{(\ell)}$ is given by
\begin{equation*}
{\mathbb G}^{(2)}_{\ell,I,\eps} (\lambda) =\sum_{k=0}^\infty
\sum\limits_{\substack{\gamma\in\FF_{I,Q}^{(\ell)} \\ q_k \leqslant
\lambda Q<q_{k+1}}} \frac{1-\eps q}{q(kq+q^\prime)} \cdot
\frac{1}{1+\gamma^2} +O\left(\frac{1}{Q}\right)  =S_{I,Q} (\lambda)
+O\left(\frac{1}{Q}\right).
\end{equation*}

When $\gamma^\prime \in \FF_{I,Q}^{(\ell)}$ similar bookkeeping with
$q_k^\prime =kq^\prime+q$, $a_k^\prime=ka^\prime +a$, $u_k
=\frac{a_k^\prime+\eps}{q_k^\prime}$, provides the contribution
\begin{equation*}
{\mathbb G}_{\ell,I,\eps}^{(3)}(\lambda)= \sum_{k=0}^\infty
\sum\limits_{\substack{\gamma^\prime\in\FF_{I,Q}^{(\ell)} \\
q^\prime_k \leqslant \lambda Q<q^\prime_{k+1}}} \left( \arctan
\gamma^\prime -\arctan u_k\right)= \sum_{k=0}^\infty
\sum\limits_{\substack{\gamma^\prime\in\FF_{I,Q}^{(\ell)} \\
q^\prime_k \leqslant \lambda Q<q^\prime_{k+1}}} \frac{1-\eps
q^\prime}{q^\prime (kq^\prime +q)}\cdot \frac{1}{1+\gamma^{\prime
2}}+O\left(\frac{1}{Q}\right).
\end{equation*}

\begin{figure}[ht]
\begin{center}
\unitlength 0.6mm
\begin{picture}(100,170)(0,0)
\texture{ccc 0}
\shade\path(22.2222,100)(33.3333,100)(33.3333,66.6666)(22.2222,100)
\texture{c 0}
\shade\path(33.3333,100)(33.3333,66.6666)(66.6666,33.3333)(33.3333,100)

\path(0,100)(100,100)(100,0)(0,100)
\path(22.2222,100)(33.3333,66.6666)
\path(33.3333,100)(66.6666,33.3333) \path(45,55)(45,100)
\dottedline{2}(66.6666,100)(66.6666,33.3333)
\dottedline{2}(45,55)(45,0)
\dottedline{2}(0,100)(0,166.6666)(22.2222,100)
\dottedline{2}(0,166.6666)(33.3333,100)
\dottedline{2}(45,77)(100,77) \dottedline{2}(45,55)(100,55)

\put(0,166.6666){\makebox(0,0){{\tiny $\bullet$}}}
\put(0,100){\makebox(0,0){{\tiny $\bullet$}}}
\put(22.2222,100){\makebox(0,0){{\tiny $\bullet$}}}
\put(33.3333,100){\makebox(0,0){{\tiny $\bullet$}}}
\put(66.6666,100){\makebox(0,0){{\tiny $\bullet$}}}
\put(100,100){\makebox(0,0){{\tiny $\bullet$}}}
\put(33.3333,66.6666){\makebox(0,0){{\tiny $\bullet$}}}
\put(66.6666,33.3333){\makebox(0,0){{\tiny $\bullet$}}}
\put(100,0){\makebox(0,0){{\tiny $\bullet$}}}
\put(45,0){\makebox(0,0){{\tiny $\bullet$}}}

\put(-9,100){\makebox(0,0){{\small $(0,1)$}}}
\put(9,166.66){\makebox(0,0){{\small $(0,\lambda)$}}}
\put(109,100){\makebox(0,0){{\small $(1,1)$}}}
\put(109,0){\makebox(0,0){{\small $(1,0)$}}}
\put(55,0){\makebox(0,0){{\small $\big( \frac{q}{Q},0\big)$}}}
\put(18,60){\makebox(0,0){{\small $\big(
\frac{\lambda-1}{k},1-\frac{\lambda-1}{k}\big)$}}}
\put(66,29){\makebox(0,0){{\small $\big(
\frac{\lambda-1}{k-1},1-\frac{\lambda-1}{k-1}\big)$}}}
\put(28,145){\makebox(0,0){{\small $y=\lambda-kx$}}}
\put(18,120){\makebox(0,0){{\small $y=\lambda-(k+1)x$}}}

\put(66,106){\makebox(0,0){{\small $\big(
\frac{\lambda-1}{k-1},1\big)$}}}
\put(37,106){\makebox(0,0){{\small $\big(
\frac{\lambda-1}{k},1\big)$}}} \put(16,106){\makebox(0,0){{\small
$\big( \frac{\lambda-1}{k+1},1\big)$}}}
\put(109,87){\makebox(0,0){{\small $J_{k,q}^{(0)}$}}}
\put(109,65){\makebox(0,0){{\small $J_{k,q}^{(1)}$}}}
\end{picture}
\end{center}
\caption{The set $\Omega_k \cap \TT$} \label{Figure6}
\end{figure}

The situation is analogous to the one encountered in \cite[Section
5]{BZ2}. Consider the ``Farey triangle" $\TT=\{(x,y)\in (0,1]^2:x+y>1\}$, the sets
\begin{equation*}
\Omega_k =\left\{ (x,y)\in\R^2: \bigg[ \frac{\lambda -y}{x}\bigg]
=k\right\} ,\quad k\geqslant 0, \qquad I_k =\left[ \frac{\lambda
-1}{k},\frac{\lambda -1}{k-1}\right) \cap [0,1),\quad k\geqslant 1,
\end{equation*}
and for $q\in QI_k$, $k\geqslant 1$, the intervals
\begin{equation*}
\begin{split}
& J_{k,q}^{(0)}=( \lambda -k\eps q ,1]=\left\{ \eps q^\prime :
(\eps q ,\eps q^\prime ) \in \Omega_{k-1}\cap \TT\right\}
\subseteq ( 1-\eps q ,1],\\
& J_{k,q}^{(1)}=( 1-\eps q,\lambda -k\eps q ]=\left\{ \eps
q^\prime : ( \eps q,\eps q^\prime ) \in \Omega_k \cap \TT\right\}
\subseteq ( 1-\eps q,1].
\end{split}
\end{equation*}
Denote
\begin{equation*}
f_{k,q}(q^\prime,a)=\frac{1-\eps q}{q(kq+q^\prime)} \cdot
\frac{1}{1+\left( \frac{a}{q}\right)^2},\quad
g_{k,q}(x,y):=f_{k,q}(x,q-y), \quad k\geqslant 0.
\end{equation*}

Assume first $\lambda \geqslant 2$. Then $\min\{ k:\Omega_k \cap
\TT\}\neq \emptyset\}=[\lambda]-1\geqslant 1$, $\min\{ k:I_k \neq
\emptyset\}=[\lambda]\geqslant 2$, and we have
\begin{equation*}
S_{I,Q}(\lambda)=\sum_{k=1}^\infty
\sum\limits_{\substack{\gamma\in\FF_{I,Q}^{(\ell)} \\
(q,q^\prime)\in Q(\Omega_k \cap \TT)}} f_{k,q}(q^\prime ,a)
=\sum_{k=2}^\infty \sum_{\substack{q\in QI_k \\ \gcd (\ell,q)=1}}
\Big( S_{I,Q,k} (\lambda,q)+T_{I,Q,k} (\lambda,q)\Big),
\end{equation*}
with
\begin{equation*}
S_{I,Q,k} (\lambda,q)=\sum\limits_{\substack{(a,q^\prime) \in
qI \times QJ_{k,q}^{(1)} \\
-aq^\prime \equiv 1\hspace{-5pt}\pmod{q},\ \ell \mid (q-a)}}
f_{k,q}(q^\prime ,a),\qquad T_{I,Q,k}(\lambda,q)=
\sum\limits_{\substack{(a,q^\prime) \in qI\times QJ_{k,q}^{(0)}\\
-aq^\prime \equiv 1\hspace{-5pt} \pmod{q},\, \ell\mid (q-a)}}
f_{k-1,q}(q^\prime ,a).
\end{equation*}
Taking $x=q^\prime$, $y=q-a=\ell z\in q(1-I)$, and $\bar{\ell}$
the multiplicative inverse of $\ell\hspace{-3pt}\pmod{q}$, we can
write
\begin{equation*}
\begin{split}
S_{I,Q,k}(\lambda,q) & =\sum\limits_{\substack{(x,y)\in\in QJ_{k,q}^{(1)}\times q(1-I),\\
xy\equiv 1 \hspace{-5pt} \pmod{q},\, \ell \mid y}} f_{k,q} (x,q-y)=
\sum\limits_{\substack{(x,z)\in QJ_{k,q}^{(1)}\times (q/\ell)
(1-I),\\ \gcd(x,q)=1,\,xz\equiv \bar{\ell}\hspace{-5pt} \pmod{q}}}
g_{k,q} (x,\ell z),\\
T_{I,Q,k}(\lambda,q) & =\sum\limits_{\substack{(x,z)\in
QJ_{k,q}^{(0)}\times (q/\ell) (1-I),\\
\gcd(x,q)=1,\,xz\equiv \bar{\ell} \hspace{-5pt} \pmod{q}}}
g_{k-1,q}(x,\ell z).
\end{split}
\end{equation*}
Since $(q,x)\in Q\Omega_k$, we have $Q-(k+1)q<x\leqslant
\lambda-kq$, so $(\lambda-1)Q \leqslant \lambda Q-q<kq+x\leqslant
t$ and one finds
\begin{equation*}
\| g_{k,q}\|_\infty \leqslant \frac{1}{qQ},\qquad \| \nabla
g_{k,q}\|_\infty \leqslant \frac{3}{q^2 Q},\qquad \forall k\geqslant 1.
\end{equation*}
Since $\ell \mid (q-a)$ and $\gcd (q,q-a)=1$ we have $\gcd
(\bar{\ell},q)=1$. The length of each of the intervals
$QJ_{k,q}^{(0)}$ and $QJ_{k,q}^{(1)}$ is less than $q$, so we can
apply Lemma \ref{L2.3} with $T=[Q^{c_1}]$ to find
\begin{equation}\label{3.3}
\begin{split}
S_{I,Q,k}(\lambda,q) & =\frac{\varphi (q)}{q^2}
\iint_{QJ_{k,q}^{(1)}\times \frac{q}{\ell} (1-I)} g_{k,q} (x,\ell
z)\, dx\, dz+\EE^{(1)}_{k,q,\ell} \\ &
=\frac{c_I}{\ell}\cdot\frac{\varphi(q)}{q^2}\cdot\frac{Q-q}{Q}\,
\int_{Q-q}^{\lambda Q-kq} \frac{dx}{x+kq}+\EE^{(1)}_{k,q,\ell},
\end{split}
\end{equation}
with
\begin{equation}\label{3.4}
\EE^{(1)}_{k,q,\ell} \ll_{\delta,\ell} T^2 \,\frac{1}{qQ}\,
q^{\frac{1}{2}+\delta} +T\, \frac{1}{q^2 Q}\,
q^{\frac{3}{2}+\delta} + \frac{q\, \frac{Q}{Q^c}\cdot\frac{1}{q^2
Q}}{T} \ll Q^{2c_1-1} q^{-\frac{1}{2}+\delta}
+\frac{Q^{-c-c_1}}{q} .
\end{equation}
A similar argument leads to
\begin{equation}\label{3.5}
T_{I,Q,k}(\lambda,q)=\frac{c_I}{\ell}\cdot
\frac{\varphi(q)}{q^2}\cdot\frac{Q-q}{Q} \int_{\lambda
Q-kq}^Q \frac{dx}{x+(k-1)q}+O_{\delta,\ell} \left( Q^{2c_1-1}
q^{-\frac{1}{2}+\delta} +\frac{Q^{-c-c_1}}{q}\right).
\end{equation}
But
\begin{equation*}
\int_{Q-q}^{\lambda Q-kq} \frac{dx}{x+kq}+\int_{\lambda Q-kq}^{Q}
\frac{dx}{x+(k-1)q}= \int_{Q+(k-1)q}^{\lambda
Q}\frac{du}{u}+\int_{\lambda Q-q}^{Q+(k-1)q} \frac{du}{u}= \ln
\frac{\lambda Q}{\lambda Q-q},
\end{equation*}
hence \eqref{3.3}--\eqref{3.5} yield
\begin{equation}\label{3.6}
S_{I,Q,k}(\lambda,q)+T_{I,Q,k}(\lambda,q)=\frac{c_I}{\ell}\cdot\frac{\varphi(q)}{q^2}
\cdot\frac{Q-q}{Q}\,\ln\frac{\lambda Q}{\lambda
Q-q}+O_{\delta,\ell} \left( Q^{2c_1-1} q^{-\frac{1}{2}+\delta}
+\frac{Q^{-c-c_1}}{q}\right).
\end{equation}
The intervals $I_k$ are disjoint, so when summing over $k$ and
$q\in QI_k$ (or in a smaller range) we are actually summing over
$q\in [1,Q]$. This way in \eqref{3.6} the error will sum up to
\begin{equation*}
O_{\delta,\ell} \left( Q^{2c_1-1} \sum_{q\leqslant Q}
q^{-\frac{1}{2}+\delta} +Q^{-c-c_1} \sum_{q\leqslant Q} \frac{1}{q}\right)
= O_{\ell,\delta} \Big( Q^{\delta-\theta (c,c_1)}\Big),
\end{equation*}
while the main term will sum up to
\begin{equation*}
M_{\ell,I}(Q)=\frac{c_I}{\ell} \sum_{\substack{
q\leqslant Q \\ \gcd (\ell,q)=1}} \frac{\varphi(q)}{q} \, W(q),
\end{equation*}
where
\begin{equation*}
W(q)=\frac{Q-q}{qQ} \, \ln \frac{\lambda Q}{\lambda Q-q},\qquad
q\in [1,Q],
\end{equation*}
with $\| W\|_\infty \ll \frac{\lambda}{Q}$ and $T_0^Q W\ll
\frac{\lambda}{Q}$. Lemma \ref{L2.1} now provides
\begin{equation*}
M_{\ell,I}(Q)=\frac{c_IC(\ell)}{\ell} \int_0^Q W(q)\,
dq+O\left(\frac{\lambda\ln Q}{Q}\right),
\end{equation*}
and we proved

\begin{proposition}\label{P3.1}
For any $\delta >0$, uniformly in $\lambda$ on compacts of $[2,\infty )$,
\begin{equation}\label{3.7}
{\mathbb G}^{(2)}_{\ell,I,\eps}(\lambda)=\frac{c_IC(\ell)}{\ell}
\int_0^1 \frac{1-u}{u}\ \ln \frac{\lambda}{\lambda -u}\
du+O_{\delta,\ell} \Big(\eps^{-\delta +\theta (c,c_1)}\Big).
\end{equation}
An identical formula holds for ${\mathbb
G}^{(3)}_{\ell,I,\eps}(\lambda)$.
\end{proposition}

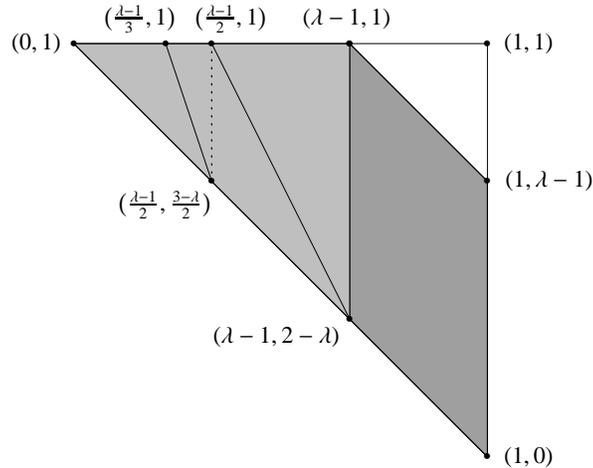
\begin{figure}[ht]
\begin{center}
\unitlength 0.55mm
\begin{picture}(100,110)(0,0)
\texture{ccc 0}
\shade\path(66.6666,100)(100,66.6666)(100,0)(66.6666,33.3333)(66.6666,100)

\texture{c 0}
\shade\path(0,100)(66.6666,33.3333)(66.6666,100)(0,100)
\path(0,100)(100,100)(100,0)(0,100)
\path(22.2222,100)(33.333,66.6666)
\path(33.333,100)(66.666,33.3333) \path(66.666,100)(100,66.6666)
\dottedline{2}(33.3333,100)(33.3333,66.6666)

\put(0,100){\makebox(0,0){{\tiny $\bullet$}}}
\put(22.2222,100){\makebox(0,0){{\tiny $\bullet$}}}
\put(33.3333,100){\makebox(0,0){{\tiny $\bullet$}}}
\put(66.6666,100){\makebox(0,0){{\tiny $\bullet$}}}
\put(100,100){\makebox(0,0){{\tiny $\bullet$}}}
\put(33.3333,66.6666){\makebox(0,0){{\tiny $\bullet$}}}
\put(66.6666,33.3333){\makebox(0,0){{\tiny $\bullet$}}}
\put(66.6666,100){\makebox(0,0){{\tiny $\bullet$}}}
\put(100,0){\makebox(0,0){{\tiny $\bullet$}}}
\put(100,66.6666){\makebox(0,0){{\tiny $\bullet$}}}

\put(-9,100){\makebox(0,0){{\small $(0,1)$}}}
\put(110,100){\makebox(0,0){{\small $(1,1)$}}}
\put(110,0){\makebox(0,0){{\small $(1,0)$}}}
\put(115,66.6666){\makebox(0,0){{\small $(1,\lambda-1)$}}}
\put(66,106){\makebox(0,0){{\small $(\lambda-1,1)$}}}
\put(49,29){\makebox(0,0){{\small $(\lambda-1,2-\lambda)$}}}
\put(22,61){\makebox(0,0){{\small $\big(
\frac{\lambda-1}{2},\frac{3-\lambda}{2}\big)$}}}
\put(38,106){\makebox(0,0){{\small $\big(
\frac{\lambda-1}{2},1\big)$}}} \put(16,106){\makebox(0,0){{\small
$\big( \frac{\lambda-1}{3},1\big)$}}}
\end{picture}
\end{center}
\caption{The set $\cup_{k=1}^\infty \Omega_k \cap \TT$ when
$1<\lambda <2$} \label{Figure7}
\end{figure}

When $1<\lambda< 2$ the contribution to $S_{I,Q}(\lambda)$ of
$\gamma$ with $q+q^\prime \leq \lambda Q$ is
\begin{equation*}
\begin{split}
A_{I,Q}(\lambda) & =\sum_{k=1}^\infty
\sum\limits_{\substack{\gamma\in \FF_{I,Q}^{(\ell)} \\ (q,q^\prime
)\in Q(\Omega_k \cap \TT)}} f_{k,q}(q^\prime,a) \\ &
=\sum_{k=2}^\infty \sum_{\substack{q\in QI_k \\ \gcd (\ell,q)=1}}
\Bigg(\sum\limits_{\substack{(a,q^\prime) \in qI\times Q J_{k,q}^{(1)}\\
aq^\prime \equiv -1\hspace{-5pt}\pmod{q},\, \ell\mid (q-a)}}
f_{k,q}(q^\prime,a)+\sum\limits_{\substack{(a,q^\prime) \in qI\times Q J_{k,q}^{(0)}\\
aq^\prime \equiv -1\hspace{-5pt}\pmod{q},\, \ell\mid (q-a)}}
f_{k-1,q}(q^\prime,a) \Bigg) \\ & \quad
+\sum_{\substack{(\lambda-1)Q\leqslant q\leqslant Q \\ \gcd
(\ell,q)=1}}
\sum\limits_{\substack{(a,q^\prime) \in qI\times Q J_{1,q}^{(1)}\\
aq^\prime \equiv -1\hspace{-5pt}\pmod{q},\, \ell\mid (q-a)}}
f_{1,q}(q^\prime,a),
\end{split}
\end{equation*}
and the contribution of $\gamma$ with $q+q^\prime > \lambda Q$ is
\begin{equation*}
B_{I,Q}(\lambda)=\sum_{\substack{(\lambda-1)Q<q\leqslant Q \\ \gcd
(\ell,q)=1}}
\sum\limits_{\substack{(a,q^\prime) \in qI\times (\lambda Q-q,Q]\\
aq^\prime \equiv -1\hspace{-5pt} \pmod{q},\,\ell \mid (q-a)}}
f_{0,q}(q^\prime ,a).
\end{equation*}
Using $\| f_{k,q}\|_\infty \leq \frac{1}{qQ}$ and $\| \nabla
f_{k,q}\|_\infty \leq\frac{3}{q^2 Q}$ on $(QJ_{k,q}^{(0)}\cup
QJ_{k+1,q}^{(1)})\times qI$ if $k\geqslant 1$, $\|
f_{0,q}\|_\infty \leqslant \frac{1}{(\lambda-1)qQ}$ and $\| \nabla
f_{0,q}\|_\infty \leqslant \frac{1}{(\lambda-1)q^2 Q}$ on
$QJ_{1,q}^{(0)} \times qI$, and summing as in the case $\lambda
>2$ we also obtain

\begin{proposition}\label{P3.2}
For every $\delta >0$,
uniformly in $\lambda$ on compacts of $(1,2]$,
\begin{equation}\label{3.8}
{\mathbb G}^{(2)}_{\ell,I,\eps}(\lambda)=\frac{c_IC(\ell)}{\ell}
\int_0^1 \frac{1-u}{u}\ \ln \frac{\lambda}{\lambda-u}\
du+O_{\delta,\ell} \Big( \eps^{-\delta +\theta (c,c_1)}\Big).
\end{equation}
An identical asymptotic formula holds for ${\mathbb
G}^{(3)}_{\ell,I,\eps}(\lambda)$.
\end{proposition}

As a result formula \eqref{3.7}
will also hold when $1<\lambda\leqslant 2$.

Consider finally the case $0<\lambda \leqslant 1$. When $q^\prime
>\lambda Q$ we have $q_k >\lambda Q$ for all $k\geqslant 0$. So the contribution of
each $\gamma\in \FF_{I,Q}^{(\ell)}$ with $q^\prime >\lambda Q$ is
in this case
\begin{equation*}
\arctan \gamma^\prime-\arctan\gamma=\frac{1}{qq^\prime}\cdot
\frac{1}{1+\gamma^2}+O\left( \frac{1}{q^2q^{\prime 2}}\right),
\end{equation*}
summing up to
\begin{equation}\label{3.9}
C_{I,Q}(\lambda) =\sum_{\substack{1\leqslant q\leqslant Q \\ \gcd
(\ell,q)=1}} \sum\limits_{\substack{(a,q^\prime) \in qI\times J_{\lambda,q}\\
\gcd(q^\prime,q)=1,\, \ell \mid (q-a) \\ aq^\prime
\equiv -1\hspace{-5pt} \pmod{q}}} \frac{1}{qq^\prime}\cdot
\frac{1}{1+\left(\frac{a}{q}\right)^2} +O\left(\frac{1}{Q}\right) .
\end{equation}
Using the same estimates as in case $\lambda>2$, \eqref{3.9} leads
to
\begin{equation}\label{3.10}
\begin{split}
C_{I,Q}(\lambda) & =\frac{c_IC(\ell)}{\ell} \int_0^1
\frac{1}{u}\,\ln \frac{1}{\max \{ 1-u,\lambda\}}\,
du+O_{\delta,\ell} \Big( Q^{\delta-\theta (c,c_1)}\Big) \\
& =\frac{c_IC(\ell)}{\ell} \Bigg( \int_0^{1-\lambda} \frac{1}{u}\
\ln \frac{1}{1-u}\ du+\ln (1-\lambda)\ln \lambda\Bigg)+
O_{\delta,\ell} \Big( Q^{\delta-\theta (c,c_1) }\Big).
\end{split}
\end{equation}
Finally the contribution of each $\gamma\in\FF_{I,Q}^{(\ell)}$
with $q^\prime \leqslant \lambda Q$ is
\begin{equation*}
\arctan t_0-\arctan\gamma=\frac{1-\eps q}{qq^\prime}\cdot
\frac{1}{1+\gamma^2}+O\left(\frac{1}{q^2q^{\prime 2}}\right),
\end{equation*}
summing up to
\begin{equation}\label{3.11}
\begin{split}
D_{I,Q}(\lambda) & = \sum_{\substack{(1-\lambda)Q<q\leqslant Q \\
\gcd (\ell,q)=1}}
\sum\limits_{\substack{(a,q^\prime) \in qI\times (Q-q,\lambda Q]\\
\gcd (q^\prime,q)=1,\, \ell \vert (q-a) \\
aq^\prime \equiv -1\hspace{-5pt} \pmod{q}}}
\frac{1-\eps q}{qq^\prime}\cdot
\frac{1}{1+\left(\frac{a}{q}\right)^2}+O\left(\frac{1}{Q}\right) \\
& =\frac{c_IC(\ell)}{\ell} \int_{1-\lambda}^1 \frac{1-u}{u}\,
\ln\frac{\lambda}{1-u}\, du+O_{\delta,\ell} \left( \lambda
Q^{-1+\delta}+Q^{\delta+\max \{ 2c_1 -\frac{1}{2},-c-c_1\}}\right).
\end{split}
\end{equation}
From \eqref{3.10} and \eqref{3.11} we infer,
uniformly in
$\lambda$ on compacts of $(0,1]$,
\begin{equation}\label{3.12}
\begin{split}
{\mathbb G}^{(2)}_{\ell,I,\eps}(\lambda) & = C_{I,Q}(\lambda)+D_{I,Q}(\lambda) \\
& =\frac{c_IC(\ell)}{\ell} \left( \int_0^1 \frac{1}{u}\ \ln
\frac{1}{1-u}\, du -\lambda\right)+O_{\delta,\ell} \Big( \lambda
Q^{-1+\delta}+Q^{\delta+\max \{ 2c_1 -\frac{1}{2},-c-c_1\}}\Big),
\end{split}
\end{equation}
whence

\begin{proposition}\label{P3.3}
For every $\delta >0$, uniformly in $\lambda$ on compacts of $(0,1]$,
\begin{equation*}
{\mathbb G}^{(2)}_{\ell,I,\eps}(\lambda) =\frac{c_IC(\ell)(\zeta
(2)-\lambda)}{\ell} +O_{\delta,\ell} \Big( \eps^{-\delta +\theta
(c,c_1)}\Big).
\end{equation*}
An identical formula holds for ${\mathbb G}^{(3)}_{\ell,I,\eps}
(\lambda)$.
\end{proposition}

\section{End of the proof of Theorem \ref{T1}}
From Corollary \ref{C2.1} and Propositions \ref{P3.1}, \ref{P3.2}, \ref{P3.3} we
gather
\begin{equation*}
\tilde{\mathbb G}_{\ell,I,\eps}(\lambda)=c_I G_\ell
(\lambda)+O_{\delta,\ell} \Big( \eps^{-\delta+\theta
(c,c_1)}\Big),
\end{equation*}
with repartition function $G_\ell$ given by
\begin{equation*}
G_\ell (\lambda) =\begin{cases} \vspace{.1cm}
1-\frac{\lambda}{\zeta (2)}+A(\ell )H_1 (\lambda) & \mbox{\rm
if $\ \lambda\in \big( 0, \frac{1}{2}\big],$} \\
\vspace{.1cm} 1-\frac{\lambda}{\zeta (2)}+A(\ell) H_2 (\lambda) &
\mbox{\rm if $\ \lambda\in\big[
\frac{1}{2}, 1\big],$} \\
\frac{2C(\ell)}{\ell}\, H_3 (\lambda) & \mbox{\rm if $\
\lambda\in [ 1,\infty),$}
\end{cases}
\end{equation*}
where
\begin{equation*}
\begin{split}
H_1 (\lambda) & =\lambda-\zeta (2)+(\ln \lambda)\ln (1-\lambda)
+2\int_\lambda^{1-\lambda} \frac{1}{u}\ \ln\frac{1}{1-u}\ du
+2\int_{1-\lambda}^1 \frac{1-u}{u}\
\ln \frac{\lambda}{1-u}\ du, \\
& =-\lambda-(\ln \lambda)\ln (1-\lambda)+\int_{1-\lambda}^1
\frac{1}{u}\ \ln \frac{1}{1-u}\ du
-\int_0^\lambda \frac{1}{u}\ \ln\frac{1}{1-u}\ du \\ & =-\lambda ,\\
H_2 (\lambda) & =\lambda-\zeta (2)+(\ln \lambda)^2
+2\int_\lambda^1 \frac{1-u}{u}\ \ln \frac{\lambda}{1-u}\ du \\
& =3\lambda-2-\zeta (2)-(\ln \lambda)^2 +2(1-\lambda) \ln \bigg(
\frac{1}{\lambda}-1\bigg) +2\int_\lambda^1 \frac{1}{u}\ \ln
\frac{1}{1-u}\ du,\\ H_3 (\lambda) & =\int_0^1 \frac{1-u}{u}\ \ln
\frac{\lambda}{\lambda-u}\ du.
\end{split}
\end{equation*}
This establishes part (i) in Theorem \ref{T1}. Using also
\begin{equation*}
H_3^\prime (\lambda)=\int_0^1 \frac{1-u}{u} \left(
\frac{1}{\lambda}-\frac{1}{\lambda-u}\right) du=-\frac{1}{\lambda}
\int_0^1 \frac{1-u}{\lambda-u}\, du =-\frac{1}{\lambda}+\left(
1-\frac{1}{\lambda}\right) \ln \left(
1-\frac{1}{\lambda}\right),\quad \lambda>1,
\end{equation*}
it follows that the density of the repartition function $G_\ell$
is
\begin{equation*}
g_\ell (\lambda)=-G_\ell^\prime (\lambda)=\begin{cases}
\frac{1}{\zeta (2)}+A(\ell) & \mbox{\rm
if $\ \lambda\in \big( 0, \frac{1}{2}\big],$} \\
\frac{1}{\zeta (2)}+A(\ell) \Big( -3+\frac{2}{\lambda} -2\big(
\frac{1}{\lambda}-1\big)\ln \big( \frac{1}{\lambda}-1\big) \Big) &
\mbox{\rm if $\ \lambda \in \big[\frac{1}{2},1\big),$} \\
\frac{2C(\ell)}{\ell}\, \Big(
\frac{1}{\lambda}-\big(1-\frac{1}{\lambda}\big) \ln \big(
1-\frac{1}{\lambda}\big)\Big) & \mbox{\rm if $\ \lambda\in (1,\infty).$}
\end{cases}
\end{equation*}

Part (ii) of Theorem \ref{T1} can now be deduced from part (i) by
a standard approximation argument \cite{BCZ,BGZ1,BGZ2,CG1} based
on $\tau_{\ell,\eps} (\omega) \approx
\frac{q_{\ell,\eps}(\omega)}{\cos \omega}$ for $\tan\omega \in I$
with $I\subseteq [0,1]$ interval of length $\vert I\vert \asymp
\eps^c$. We skip the detailed proof, which can be easily
reconstructed from some of the arguments in the next section.

\section{The conversion from a honeycomb to a square lattice with congruence constraints}

\begin{figure}[ht]
\includegraphics*[scale=0.93, bb=90 0 280 105]{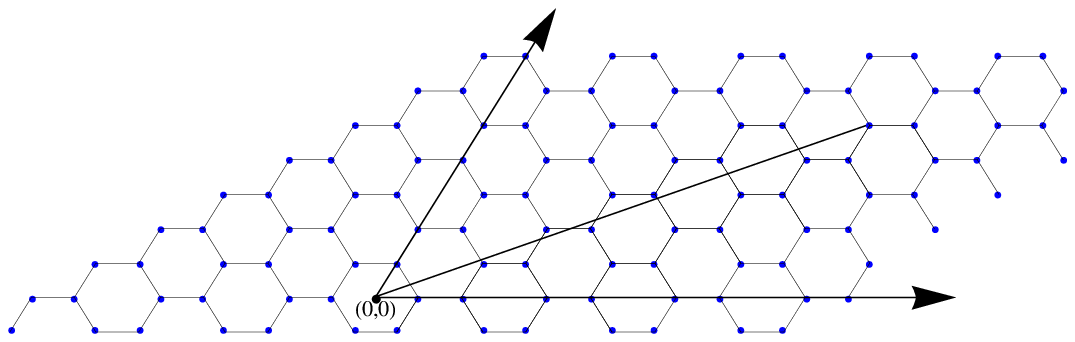}
\includegraphics*[scale=0.7, bb=-30 0 300 140]{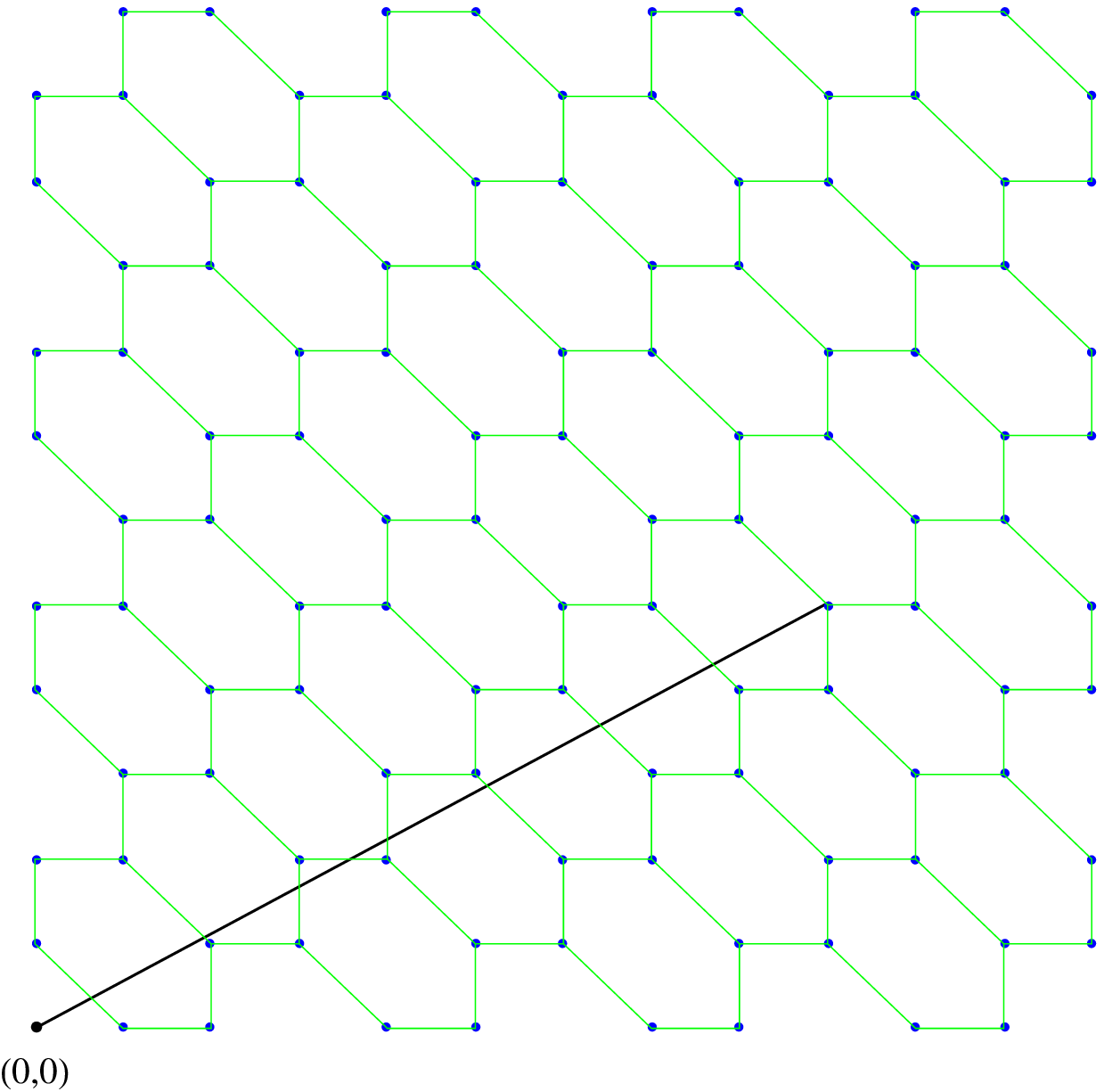}
\caption{\it \small The free path length in the honeycomb and in
the deformed honeycomb } \label{Figure8}
\end{figure}

In the situation of the honeycomb it suffices to consider $\omega
\in [0,\frac{\pi}{6}]$. The linear transformation
\begin{equation*}
\R^2 \ni (x,y)\ \ \stackrel{T}{\longrightarrow}\ \
(x^\prime,y^\prime)=\left(
x-\frac{y}{\sqrt{3}},\frac{2y}{\sqrt{3}}\right)\in \R^2
\end{equation*}
maps the first sextant
$\Gamma_+:=\{(q+\frac{a}{2},\frac{a\sqrt{3}}{2}): q,a\in \Z,
0\leqslant a\leqslant q\}$ of the grid of equilateral triangles of
side $1$ onto the first quadrant of the square lattice
$\Z^2=\{(q,a) :a,q\in\Z\}$. Elements of the subset
$\Gamma_+^{\mathrm{hex}}\subseteq \Gamma_+$ of vertices from the
honeycomb grid map to integer lattice points $(q,a)$ with
$q\nequiv a\hspace{-3pt} \pmod{3}$ (see Figure \ref{Figure8}).

\begin{figure}[ht]
\includegraphics*[scale=1.75, bb=0 0 230 70]{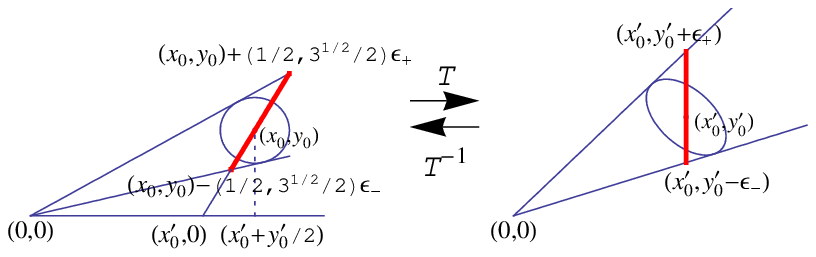}
\caption{\it \small Change of scatterers under the linear
transformation $T$ } \label{Figure9}
\end{figure}

$T$ also maps circular scatterers
$(q+\frac{a}{2},\frac{a\sqrt{3}}{2})+\eps (\cos\theta,\sin\theta)$
centered at $(x_0,y_0)=(q+\frac{a}{2},\frac{a\sqrt{3}}{2})$ to
ellipsoidal scatterers $(q,a)+\eps (\cos\theta
-\frac{\sin\theta}{\sqrt{3}},\frac{2\sin\theta}{\sqrt{3}})$
centered at $(x_0^\prime,y_0^\prime)=(q,a)$. Denote
$\omega_0=\arctan \frac{y_0}{x_0}$ and $\omega_0^\prime=\arctan
\frac{y_0^\prime}{x_0^\prime}$. Denote also $S_{\eps}=\{ \delta
(\cos\frac{\pi}{3},\sin\frac{\pi}{3}): \vert \delta\vert \leqslant
\eps\}$.

In the honeycomb with scatterers
$\{(x_0,y_0)+S_{\eps}:(x_0,y_0)\in\Gamma_+^{\,\operatorname{hex}}\}$
of same length, the "local" repartition function
\begin{equation}\label{4.1}
\tilde{\mathbb G}^{\,\operatorname{hex}}_{I,\eps} (\lambda)=\bigg|
\bigg\{ \omega\in \arctan I: \tilde{q}^{\,\operatorname{hex}}_\eps
(\omega)>\frac{\lambda}{\eps}\bigg\} \bigg| ,\qquad I\subseteq
\bigg[0,\frac{\pi}{6}\bigg]\ \ \mbox{\rm interval,}
\end{equation}
of the "horizontal" free path length
\begin{equation*}
\tilde{q}_\eps^{\,\operatorname{hex}} (\omega) =\inf \left\{ q\in
\N: \exists a\in \N,\exists (x_0,y_0)\in
\Gamma_+^{\,\operatorname{hex}} , \bigg(
\frac{q+a/2}{\cos\omega},\frac{a\sqrt{3}/2}{\sin\omega}\bigg) \in
(x_0,y_0)+S_{\eps}\right\},
\end{equation*}
turns out to be closely related with ${\mathbb
G}_{3,I,\eps}(\lambda)$, being estimated through the same
approximation procedure as for the later when $I\subseteq [0,1]$
is a short interval of length $\vert I\vert\asymp \eps^c$,
$0<\eps<1$. Indeed, the equality
\begin{equation*}
T^{-1} (x_0^\prime ,y_0^\prime+\delta)=\Bigg( x_0^\prime
+\frac{y_0^\prime+\delta}{2},\frac{(y_0^\prime
+\delta)\sqrt{3}}{2}\Bigg),
\end{equation*}
shows that $T$ maps the oblique scatterer $(x_0,y_0)+S_{\eps}$
from the honeycomb onto the vertical scatterer
$(x_0^\prime,y_0^\prime)+V_{\eps}$ from the square lattice with
$x_0^\prime \equiv y_0^\prime \hspace{-3pt} \pmod{3}$, and the line
through the origin with slope $\frac{y_0^\prime+\delta}{x_0^\prime}$
onto the line through the origin with slope $\Phi \big(
\frac{y_0^\prime+\delta}{x_0^\prime}\big)$,
where $\Phi$ is the bijection
\begin{equation*}
\Phi :[0,1]\rightarrow \Bigg[ 0,\frac{\sqrt{3}}{3}\Bigg],\qquad
\Phi (\mu)=\frac{\mu\sqrt{3}}{2+\mu},\qquad \Phi^{-1}
(x)=\frac{2x}{\sqrt{3}-x}.
\end{equation*}
In the process the condition $\tilde{q}_\eps^{\,
\operatorname{hex}} (\omega) >\frac{\lambda}{\eps}$ above is being
replaced by $q_{3,\eps} (\omega^\prime)>\frac{\lambda}{\eps}$ in
the square lattice. The only difference arises from replacing
expressions (with $x_0^\prime =q$, $y_0^\prime =a$ and $0\leqslant
\delta \leqslant \frac{1}{q}$)
\begin{equation*}
\arctan \frac{y_0^\prime+\delta}{x_0^\prime}-\arctan
\frac{y_0^\prime}{x_0^\prime}=\frac{\delta}{x_0^\prime}\cdot
\frac{1}{1+\frac{y_0^{\prime 2}}{x_0^{\prime 2}}}+O\Bigg(
\frac{\delta^2}{x_0^{\prime 2}}\Bigg)=\frac{\delta q}{q^2+a^2}
+O\Bigg( \frac{\delta^2}{q^2}\Bigg)
\end{equation*}
that collect the contribution of angles $\omega$ for ${\mathbb
G}_{3,I,\eps}(\lambda)$, by
\begin{equation*}
\begin{split}
\arctan \Phi \Bigg(
\frac{y_0^\prime+\delta}{x_0^\prime}\Bigg)-\arctan \Phi \Bigg(
\frac{y_0^\prime}{x_0^\prime}\Bigg) & =\arctan
\frac{(y_0^\prime+\delta)\sqrt{3}}{2x_0^\prime +y_0^\prime
+\delta}-\arctan \frac{y_0^\prime \sqrt{3}}{2x_0^\prime
+y_0^\prime} \\ & =\Bigg(
\frac{(y_0^\prime+\delta)\sqrt{3}}{2x_0^\prime +y_0^\prime
+\delta}-\frac{y_0^\prime \sqrt{3}}{2x_0^\prime +y_0^\prime}
\Bigg)\frac{1}{1+\Big( \frac{y_0^\prime \sqrt{3}}{2x_0^\prime
+y_0^\prime}\Big)^2} \\ & =\frac{2x_0^\prime
\delta\sqrt{3}}{(2x_0^\prime +y_0^\prime)^2} \Bigg( 1+O\bigg(
\frac{\delta}{2x_0^\prime +y_0^\prime}\bigg) \Bigg)
\frac{1}{1+\Big( \frac{y_0^\prime \sqrt{3}}{2x_0^\prime
+y_0^\prime}\Big)^2} \\ & =\frac{\delta
q\sqrt{3}}{2(q^2+aq+a^2)}+O\Bigg( \frac{\delta^2}{q^2}\Bigg).
\end{split}
\end{equation*}
The effect will only be on the main term, where $c_J=\int_J
\frac{dx}{1+x^2}$, $J\subseteq [0,1]$, will be replaced by
$c_I^{\operatorname{hex}} :=\int_{\Phi^{-1} (I)}
\frac{\sqrt{3}}{2(1+x+x^2)}\, dx$, $I\subseteq \big[
0,\frac{1}{\sqrt{3}}\big]$. In this way we obtain, uniformly in
$\lambda$ on compacts in $\R_+\setminus \{ 1\}$,
\begin{equation}\label{4.2}
\tilde{\mathbb
G}^{\,\operatorname{hex}}_{I,\eps}(\lambda)=c_I^{\operatorname{hex}}
G_3 (\lambda)+O_{\delta} \Big(\eps^{-\delta+\theta (c,c_1)}\Big).
\end{equation}
The change of variable $x=\Phi (\mu)$ gives
\begin{equation}\label{4.3}
c^\Delta_{[\tan\omega_0,\tan\omega_1]}=\frac{\sqrt{3}}{2}
\int_{\Phi^{-1} (\tan\omega_0)}^{\Phi^{-1} (\tan\omega_1)}
\frac{dx}{x^2+x+1} = \int_{\tan\omega_0}^{\tan\omega_1}
\frac{d\mu}{\mu^2+1}=\omega_1-\omega_0.
\end{equation}
In this case $\tilde{q}_\eps^{\,\operatorname{hex}} (\omega)$ and
the free path length $\tilde{\tau}_\eps^{\,\operatorname{hex}}
(\omega)$ are related (by the rule of Sines) by
\begin{equation*}
\tilde{\tau}_\eps^{\,\operatorname{hex}}
(\omega)=\frac{\sin\frac{2\pi}{3}}{\sin \big(
\frac{\pi}{3}-\omega\big)} \ \tilde{q}_\eps^{\,\operatorname{hex}}
(\omega) =\frac{\sqrt{3}}{2} \cdot
\frac{\tilde{q}_\eps^{\,\operatorname{hex}} (\omega)}{\cos
\big(\frac{\pi}{6}+\omega\big)},
\end{equation*}
so that
\begin{equation}\label{4.4}
\tilde{\mathbb P}^{\,\operatorname{hex}}_{I,\eps} (\lambda)
:=\bigg| \bigg\{ \omega\in\arctan I:
\tilde{\tau}_\eps^{\,\operatorname{hex}}
(\omega)>\frac{\lambda}{\eps}\bigg\} \bigg|  =\bigg| \bigg\{
\omega\in\arctan I: \tilde{q}^{\,\operatorname{hex}}_\eps
(\omega)>\frac{2\lambda\cos\big(
\frac{\pi}{6}+\omega\big)}{\eps\sqrt{3}} \bigg\} \bigg|.
\end{equation}
Fix $\omega_I\in\arctan I$. Using $\big|
\cos\big(\frac{\pi}{6}+\omega\big) -\cos \big(
\frac{\pi}{6}+\omega_I\big)\big| \leqslant \vert\omega
-\omega_I\vert \ll \eps^c$, \eqref{2.2}, and
$c_I^{\operatorname{hex}} \asymp \eps^c$, equalities \eqref{4.2}
and \eqref{4.4} yield
\begin{equation}\label{4.5}
\tilde{\mathbb P}^{\,\operatorname{hex}}_{I,\eps}(\lambda)
=c_I^{\operatorname{hex}} G_3 \Bigg(
\frac{2\lambda\cos\big(\frac{\pi}{6}+\omega_I)}{\sqrt{3}} \Bigg)
+O_\delta \Big( \eps^{2c}+\eps^{-\delta+\theta (c,c_1)}\Big).
\end{equation}

Let $\eps_\pm=\eps_\pm (\omega_0,\eps)$ as in Figure
\ref{Figure9}. Consider the case of scatterers
$(x_0,y_0)+S_{\omega_0,\eps}$,
$(x_0,y_0)\in\Gamma_+^{\,\operatorname{hex}}$, where
$S_{\omega_0,\eps}$ denotes the segment $\{ \delta
(\cos\frac{\pi}{3},\sin\frac{\pi}{3}):-\eps_- \leqslant \delta
\leqslant \eps_+\}$. Let
$\tilde{\tilde{\tau}}^{\,\operatorname{hex}}_\eps (\omega)$ denote
the free path length and $\tilde{\tilde{{\mathbb
P}}}^{\,\operatorname{hex}}_{I,\eps} (\lambda)=\vert \{
\omega\in\arctan I:
\tilde{\tilde{\tau}}^{\,\operatorname{hex}}_\eps
(\omega)>\frac{\lambda}{\eps} \} \vert$. Since $\vert
\cos\omega_\pm -\cos \omega_I\vert \ll \eps^c$ and $\cos
(\frac{\pi}{6}+\omega_I)>\frac{1}{2}$ there exists $C_1 >0$ such
that
\begin{equation}\label{4.6}
\eps^\prime :=\eps \left( \frac{1}{\cos (\frac{\pi}{6}+\omega_I)}
-C_1 \eps^c \right) \leqslant \eps_- \leqslant \eps_+ \leqslant
\eps^{\prime \prime} :=\eps \left( \frac{1}{\cos
(\frac{\pi}{6}+\omega_I)}+C_1 \eps^c\right).
\end{equation}
The inequalities $\tilde{\tau}^\Delta_{\eps^{\prime\prime}}
(\omega)\leqslant \tilde{\tilde{\tau}}^\Delta_\eps
(\omega)\leqslant \tilde{\tau}^\Delta_{\eps^\prime} (\omega)$ and
\eqref{4.6} combined with formula \eqref{4.5} and \eqref{2.2} lead
to
\begin{equation*}
\begin{split}
\tilde{\tilde{{\mathbb P}}}_{I,\eps}^{\,\operatorname{hex}}
(\lambda) & \leqslant \bigg| \bigg\{ \omega\in \arctan I:
\tilde{\tau}^{\,\operatorname{hex}}_{\eps^\prime}
(\omega)>\frac{\lambda}{\eps}\bigg\} \bigg| \\ & =\Bigg| \Bigg\{
\omega \in\arctan I:
\tilde{\tau}^{\,\operatorname{hex}}_{\eps^\prime} (\omega) >
\frac{\lambda}{\eps^\prime \big( \cos
(\frac{\pi}{6}+\omega_I)-C_1\eps^c\big)}\Bigg\} \Bigg|
\\ & =\tilde{{\mathbb P}}^{\,\operatorname{hex}}_\eps \left( \frac{\lambda}{\cos
(\frac{\pi}{6}+\omega_I)-C_1 \eps^c} \right) \\ &
=c_I^{\operatorname{hex}} G_3 \left( \frac{2\lambda}{\sqrt{3}}
\cdot \frac{\cos (\frac{\pi}{6}+\omega_I)}{\cos
(\frac{\pi}{6}+\omega_I)-C_1 \eps^c}\right)+O_\delta
\Big(\eps^{2c}+\eps^{-\delta+\theta
(c,c_1)}\Big) \\
& =c_I^{\operatorname{hex}} G_3 \left(
\frac{2\lambda}{\sqrt{3}}\right)+O_\delta \Big(
\eps^{2c}+\eps^{-\delta+\theta (c,c_1)}\Big),
\end{split}
\end{equation*}
and to a similar lower bound for $\tilde{\tilde{{\mathbb
P}}}_{I,\eps}^{\,\operatorname{hex}} (\lambda)$, and so we get
\begin{equation}\label{4.7}
\tilde{\tilde{{\mathbb P}}}_{I,\eps}^{\,\operatorname{hex}}
(\lambda) =c_I^{\operatorname{hex}} G_3 \left(
\frac{2\lambda}{\sqrt{3}}\right)+O_\delta \Big(
\eps^{2c}+\eps^{-\delta+\theta (c,c_1)}\Big).
\end{equation}

The trivial inequality $\vert \tau_\eps^\Delta
(\omega)-\tilde{\tilde{\tau}}_\eps^\Delta (\omega)\vert \leqslant
2\eps$ and \eqref{4.7} now provide the formula
\begin{equation}\label{4.8}
{\mathbb P}^{\,\operatorname{hex}}_{I,\eps} (\lambda)
=c_I^{\operatorname{hex}} G_3 \left(
\frac{2\lambda}{\sqrt{3}}\right)+O_\delta \Big(
\eps^{2c}+\eps^{-\delta+\theta (c,c_1)}\Big)
\end{equation}
for the repartition function ${\mathbb
P}^{\,\operatorname{hex}}_{I,\eps} (\lambda)=\vert\{
\omega\in\arctan I: \tau_\eps^{\,\operatorname{hex}}
(\omega)>\frac{\lambda}{\eps}\}\vert$ of $\eps
\tau_\eps^{\,\operatorname{hex}}$.

Finally we choose a partition $\big[
0,\frac{1}{\sqrt{3}}\big]=\cup_{j=1}^N I_j$ with intervals $I_j$
of equal size $\frac{1}{N}\asymp \eps^c$. Applying \eqref{4.8} to
each individual interval $I_j$ with $c=c_1=\frac{1}{8}$ and
summing over $j$ we find
\begin{equation}\label{4.9}
{\mathbb P}^{\,\operatorname{hex}}_{[0,1/\sqrt{3}],\eps}
(\lambda)=\sum_{j=1}^N c_{I_j}^{\,\operatorname{hex}} G_3 \left(
\frac{2\lambda}{\sqrt{3}}\right) +O_\delta \Big(
\eps^{\frac{1}{8}-\delta}\Big)=\frac{\pi}{6}G_3 \left(
\frac{2\lambda}{\sqrt{3}}\right) +O_\delta \Big(
\eps^{\frac{1}{8}-\delta}\Big).
\end{equation}
Theorem \ref{T2} now follows immediately from \eqref{4.9} and
obvious symmetry properties of the honeycomb.

\section*{Acknowledgments}
We would like to thank the referee for careful reading and useful comments.
The second author thanks Department of Mathematics, UIUC, for hospitality
during his Fall 2007 visit, when most of this research has been completed.

\end{document}